\documentclass[preprint, 11pt]{elsarticle}

\usepackage[T1]{fontenc}
\usepackage[utf8]{inputenc}
\usepackage[hmargin=0.9in]{geometry}
\usepackage[babel,german=quotes]{csquotes}
\usepackage{subcaption}
\usepackage{graphicx}
\usepackage[colorlinks=true,citecolor=blue,urlcolor=blue]{hyperref}
\usepackage{caption}
\usepackage{amsmath}          
\usepackage{amsthm}           
\usepackage{amssymb}          
\usepackage{bm, amsfonts}     
\usepackage{xcolor}
\usepackage{fixmath}
\usepackage{tikz}
\usepackage{bm}
\usepackage{makecell}

\newcommand{\pd}[2]{\frac{\partial #1}{\partial #2}}

\renewcommand{\vec}[1]{\ensuremath{\mathbf{#1}}}

\newcommand{\dx}{\mathrm{d}\mathbf{x}}
\newcommand{\ds}{\mathrm{d}\mathbf{s}}

\newcommand{\bx}{\vec{x}}

\renewcommand{\iota}{\mathrm{i}}

\newcommand{\R}{\mathbb{R}}

\newcommand{\beq}{\begin{equation}}
\newcommand{\eeq}{\end{equation}}

\newcommand{\mina}[1]{\min\left(#1\right)}

\newcommand{\mins}[1]{\min(#1)}





\newcommand{\doi}[1]{\textrm{\textsc{doi:}} \href{http://dx.doi.org/#1}{\nolinkurl{#1}}}

\newdefinition{remark}{Remark}

\makeatletter
\def\ps@pprintTitle{%
  \let\@oddhead\@empty
  \let\@evenhead\@empty
  \def\@oddfoot{
    \footnotesize\itshape
    \ifx\@journal\@empty Elsevier
    \else\@journal\fi
    \hfill\today
  }%
  \let\@evenfoot\@oddfoot}
\makeatother

\begin{document}

\begin{frontmatter} 
  \title{Dissipation-based WENO stabilization of high-order finite element methods for scalar conservation laws}

\author{Dmitri Kuzmin\corref{cor1}}
\ead{kuzmin@math.uni-dortmund.de}

\cortext[cor1]{Corresponding author}

\author{Joshua Vedral}
\ead{joshua.vedral@math.tu-dortmund.de}

\address{Institute of Applied Mathematics (LS III), TU Dortmund University\\ Vogelpothsweg 87,
  D-44227 Dortmund, Germany}

\journal{Submitted to J.~Comput.~Phys.~special issue in honor of Prof.~Roland Glowinski}

\begin{abstract}
  We present a new perspective on the use of weighted essentially nonoscillatory (WENO) reconstructions in high-order methods for scalar hyperbolic conservation laws.
  The main focus of this work is on nonlinear stabilization of
  continuous Galerkin (CG) approximations. The
  proposed methodology also provides an interesting alternative to
  WENO-based limiters for  discontinuous Galerkin (DG) methods.
  Unlike Runge--Kutta DG schemes that
overwrite finite element solutions with WENO reconstructions, our approach uses a reconstruction-based smoothness
sensor to blend the numerical viscosity
operators of high- and low-order stabilization terms.
The so-defined  WENO approximation introduces low-order nonlinear
diffusion in the vicinity of shocks, while preserving the high-order accuracy
of a linearly stable baseline discretization in regions
where the exact solution is sufficiently smooth. The underlying reconstruction
procedure performs Hermite interpolation on stencils consisting of a
mesh cell and its neighbors. The
amount of numerical dissipation depends on the relative differences between partial
derivatives of reconstructed 
candidate polynomials and those of the underlying finite element approximation.
All derivatives are taken
into account by the employed smoothness sensor. To assess the accuracy
of our CG-WENO scheme, we derive error estimates and perform
numerical experiments. In particular, we prove that the consistency
error of the nonlinear stabilization is of the order $p+1/2$,
where $p$ is the polynomial degree. This estimate is optimal for general
meshes. For uniform meshes and smooth
exact solutions, the experimentally observed rate of convergence is as high
as $p+1$.

\end{abstract}
\end{frontmatter}

\section{Introduction}

Many reconstruction-based finite volume and discontinuous Galerkin (DG) methods for hyperbolic problems are designed to adaptively select or blend polynomial approximations corresponding to alternative stencils. 
Classical representatives of such approaches, such as the \emph{essentially nonoscillatory} (ENO) schemes developed by Harten et al. \cite{harten1987}, use the resulting reconstructions to calculate the Riemann data for high-order extensions of Godunov's method. Typical requirements for the adaptive selection of stencils/weights include crisp resolution of discontinuities, lack of spurious ripples, and uniformly high accuracy for reconstructions of smooth functions. In weighted ENO (WENO) schemes, convex combinations of candidate polynomials are constructed using normalized smoothness indicators as nonlinear weights \cite{shu1998,shu2009}.  
The WENO methodology was introduced by Liu et al.~\cite{liu1994} in the context of one-dimensional finite volume approximations. The highly cited paper by Jiang and Shu \cite{jiang1996} addressed the aspects of efficient implementation on structured grids in 1D and 2D. The fifth-order WENO scheme proposed in \cite{jiang1996} and the underlying smoothness indicator have greatly influenced
further research efforts in the field. In contrast to ENO algorithms that use binary weights \mbox{(0 or 1)}, the WENO approach produces numerical fluxes that depend continuously on the data.

The development of modern ENO/WENO reconstruction tools for unstructured grids
was initiated by Abgrall \cite{abgrall1994} and  Friedrich \cite{friedrich1998}.
Qiu and Shu \cite{qiu2005} found that local postprocessing based on WENO 
reconstructions is an excellent alternative to traditional slope limiting in
Runge--Kutta discontinuous Galerkin (DG) methods. Tailor-made algorithms
for Hermite WENO (HWENO) limiting in the DG setting
were proposed in \cite{luo2007,luo2012,zhu2009,zhu2016,zhu2017}. 
Zhang and Shu \cite{zhang2010b,zhang2010c,zhang2011,zhang2012} constrained
their WENO reconstructions using a simple scaling
limiter that enforces appropriate maximum
principles at Legendre Gauss--Lobatto quadrature points and ensures
positivity preservation for the evolved cell averages. Detailed reviews of the
literature on WENO-DG schemes can be found in \cite{zhang2011,shu2016}.

 Finite element methods based on high-order continuous Galerkin (CG)
 approximations do not support the possibility of directly adjusting
 the gradients of the approximate solution in troubled cells. However,
 spurious oscillations and violations of discrete maximum principles
 (DMPs) can be avoided by adding stabilization terms to the baseline
 discretization. 
 The framework of \emph{algebraic flux correction}
 (AFC) makes it possible to guarantee preservation of global and/or
 local bounds for scalar quantities of interest. The AFC schemes
 reviewed in \cite{kuzmin2012a,kuzmin2023} use artificial diffusion
 operators and flux limiting techniques for this purpose.
 A potential disadvantage compared to WENO-DG approaches
 lies in the fact that local DMPs impose a
 second-order accuracy barrier \cite{zhang2011}, while limiters based on global
 DMP constraints may fail to prevent nonphysical solution behavior.
 A possible remedy is the use of smoothness indicators to blend
 local and global bounds (as proposed, e.g., in \cite{hajduk2020b}).
 
 In this paper, we use HWENO reconstructions to design smoothness
 sensors for algebraic stabilization of high-order CG schemes. The
 proposed method introduces high-order dissipation in smooth regions
 and low-order dissipation in the neighborhood of discontinuities.
 This design philosophy traces its origins to the classical
 Jameson--Schmidt--Turkel (JST) scheme \cite{jameson2017,jameson1981}.
 The JST smoothness indicator for 1D schemes 
 is defined using finite difference approximations
 to the second and first derivatives. 
 Extensions to CG-$\mathbb{P}_1$ approximations
 on unstructured meshes can be found, e.g., in
 \cite{badia2014,barrenechea2017a,selmin1993}.
 Barrenechea et al. \cite{barrenechea2017c} perform in-depth theoretical
 investigations of a first-order artificial diffusion method combined 
 with a second-order local projection stabilization
 (LPS) scheme.
 
 The construction of
 JST-type stabilization terms for CG and DG methods using
 polynomials of degree $p>1$ is a more delicate issue.
 In general, it is essential to ensure that
 \begin{itemize}
\item[(i)] high-order (HO) stabilization does not degrade the rates
 of convergence to smooth solutions and vanishes
 in the neighborhood of discontinuities;
\item[(ii)] low-order (LO) dissipation is strong enough to
  suppress spurious
  oscillations but vanishes if the
  exact solution is a polynomial of degree $p$.
 \end{itemize}
 Upper bounds for the viscosity parameters of HO and LO
 stabilization operators can be obtained by estimating the
 maximum wave speed as in
\cite{abgrall2006,kuzmin2020g,lohmann2017}.
Multiplication of the two components by nonnegative weights that
add up to one makes it possible to select any convex combination
of HO and LO terms. In our method, the HO weight
$\gamma_e\in[0,1]$ of cell $e$ depends on the difference between (the
derivatives of) the evolved
finite element solution and a WENO reconstruction.
The LO weight is given by $1-\gamma_e$. The requirements
(i) and (ii) imply that we should use $\gamma_e=0$ around
shocks and $\gamma_e=1$ in cells belonging to smooth regions.
Adopting these general design principles, we construct and analyze
a nonlinear blend of HO and LO stabilization terms. The
results of numerical experiments for standard 1D and 2D test
problems illustrate the excellent shock-capturing capabilities of our
dissipation-based CG-WENO scheme. Optimal convergence rates
are attained for smooth data.

In the next section, we present the generic form of a
stabilized CG method for a scalar conservation law.
Section~3 introduces a nonlinear blend of 
HO and LO stabilization terms. In Section 4, we define
a WENO-based smoothness indicator. Some details of
the employed reconstruction procedures are given in Section~5.
The analysis presented in Section~6 yields an optimal
$\mathcal O(h^{p+1/2})$ estimate of the consistency
error.
In the last two sections, we show numerical examples
and draw conclusions.

\section{Stabilized Galerkin discretizations}

Let $u(\mathbf{x},t)$ be a scalar conserved quantity depending on
the space location $\mathbf{x}\in\bar\Omega$ and time instant $t\ge 0$.
The Lipschitz boundary of the spatial domain $\Omega\subset\R^d,\ d\in\{1,2,3\}$
is denoted by~$\Gamma=\partial\Omega$.
Imposing periodic boundary conditions on $\Gamma$, we consider
the initial value problem
\begin{subequations}\label{ibvp}
\begin{alignat}{3}
  \pd{u}{t}+\nabla\cdot\vec{f}(\mathbf{x},u)&=0 &&\quad\mbox{in}\ \Omega\times (0,T),
    \label{ibvp-pde}\\
     u(\cdot,0)&=u_0 &&\quad\mbox{in}\ \Omega,\label{ibvp-ic}
\end{alignat}
\end{subequations}
where $u_0$ is the initial data and $\mathbf{f}(\mathbf{x},u)$ is
the flux function of the conservation law. For example, if $u$ is
advected by a given velocity field $\mathbf{v}=\mathbf{v}(\mathbf{x})$,
then $\vec{f}(\mathbf{x},u)=\mathbf{v}(\mathbf{x})u$. In other
hyperbolic problems of the form \eqref{ibvp}, the flux vector $\mathbf{f}$
may be independent of $\mathbf{x}$ but depend nonlinearly on $u$.

We discretize \eqref{ibvp} in space using the continuous
Galerkin method on a conforming affine mesh $\mathcal T_h=\{
K_{1},\ldots,K_{E_h}\}$. For simplicity, we assume that
$\bigcup_{K\in\mathcal T_h}=\bar\Omega$. The mesh size
corresponding to $\mathcal T_h$ is defined by
$h=\max_{K\in\mathcal T_h}h_K$, where $h_K=\mathrm{diam}(K)$.
We seek an approximate solution
$$
u_h=\sum_{j=1}^{N_h}u_j\varphi_j
$$
in a finite element space $V_h$ spanned by Lagrange, Bernstein, or
Legendre Gauss-Lobatto (LGL)
basis functions $\varphi_1,\ldots,\varphi_{N_h}$. The methodology
to be presented below is independent of the basis. The polynomial
degree of the local finite element approximation $u_h^e=u_h|_{K_e}$ is
denoted by $p$.

The standard Galerkin discretization of \eqref{ibvp} leads to
the semi-discrete problem
\beq\label{semi-Gal}
\sum_{e=1}^{E_h}\int_{K_e}w_h\left(\pd{u_h}{t}
+\nabla\cdot\vec{f}(\mathbf{x},u_h)\right)\dx=0\qquad\forall w_h\in V_h.
\eeq
It is well known that this spatial semi-discretization may not
exhibit optimal convergence behavior even if the exact solution $u$
is smooth. Quarteroni
and Valli \cite[14.3.1]{quarteroni1994} prove that
$\|u-u_h\|_{L^2(\Omega)}=\mathcal O(h^p)$
for linear advection problems and general meshes. For a properly
stabilized CG method, the $L^2$ error is
$\mathcal O(h^{p+1/2})$; see, e.g.,
 \cite{burman2010,hartmann2008}.  All schemes that we consider
below can be written as
\beq\label{semi-stab}
\sum_{e=1}^{E_h}\int_{K_e}w_h\left(\pd{u_h}{t}
+\nabla\cdot\vec{f}(\mathbf{x},u_h)\right)\dx+
\sum_{e=1}^{E_h}s_h^e(u_h,w_h)=0\qquad\forall w_h\in V_h.
\eeq
We discuss some old and new definitions of the local stabilization
operator $s_h^e(\cdot,\cdot)$ in the next section.

Discretization in time can be performed, e.g., using a strong
stability preserving (SSP) Runge--Kutta method \cite{gottlieb2001}.
The $L^2$ error analysis for fully discrete problems confirms the
need for adding stabilization terms or using dissipative time
stepping (as in Taylor--Galerkin methods \cite{donea2003}).

\section{Dissipation-based stabilization}

The simplest way to stabilize a scheme that produces 
spurious oscillations is to add large
amounts of isotropic artificial diffusion. We define the
corresponding low-order stabilization operator
\beq\label{stab:LO}
s_h^{e,L}(u_h,w_h)=\nu_e\int_{K_e}\nabla w_h\cdot \nabla u_h\dx
\eeq
using the viscosity parameter
\[
\nu_e=\frac{\lambda_e h_e}{2p},
\]
where $h_e$ is the local mesh size and $\lambda_e=\|\mathbf{f}'(u_h)\|_{L^\infty(K_e)}$ is an upper bound for the local wave speed.

The use of \eqref{stab:LO} may be appropriate for cells located in steep
front regions. However, stabilization of
this kind introduces an $\mathcal O(h^{1/2})$ consistency error (see Section 6).
A modified version \cite{kuzmin2020g,lohmann2017} of the two-level
variational multiscale (VMS) method proposed by John et al. \cite{john2006}
replaces \eqref{stab:LO} with 
\beq\label{stab:HO}
s_h^{e,H}(u_h,w_h)=\nu_e
\int_{K_e}\nabla w_h\cdot(\nabla u_h-\mathbf{g}_h)\dx,
\eeq
where $\mathbf{g}_h$ is a continuous approximation to $\nabla u_h$.
 In essence, this stabilization technique adds a linear
 antidiffusive correction to \eqref{stab:LO}. If
 $\|\mathbf{g}_h-\nabla u_h\|_{L^2(\Omega)}=\mathcal O(h^p)$,
 then the consistency error is $\mathcal O(h^{p+1/2})$, as
 we show in Section 6 for a symmetric counterpart of this stabilization operator.
 
 Lohmann et al.
 \cite{lohmann2017} discovered an interesting
 relationship of \eqref{stab:HO} to a consistent
 streamline upwind Petrov--Galerkin (SUPG) method
 \cite{brooks1982,burman2010}. It turned out that the
 two approaches are equivalent in 1D if
 $\mathbf{g}_h\in (V_h)^d$ is defined as the
consistent-mass $L^2$ projection of $\nabla u_h$, i.e., if
\beq\label{gradrec:mcons}
\sum_{e=1}^{E_h}
 \int_{K_e} w_h(\mathbf{g}_h-\nabla u_h)\dx = 0 \qquad \forall w_h\in V_h.
 \eeq
 Substituting $w_h\in\{\varphi_1,\ldots,\varphi_{N_h}\}$, we find that
  the nodal values of 
 $\mathbf{g}_h=\sum_{j=1}^{N_h}\mathbf{g}_j\varphi_j$ satisfy
 \[
 \sum_{j=1}^{N_h}m_{ij}\mathbf{g}_j=\sum_{j=1}^{N_h}\mathbf{c}_{ij}u_j,
 \qquad i=1,\ldots,N_h,
 \]
 where
 \[
m_{ij}=\sum_{e=1}^{E_h}
 \int_{K_e}\varphi_i\varphi_j\dx,\qquad
 \mathbf{c}_{ij}=\sum_{e=1}^{E_h}
 \int_{K_e}\varphi_i\nabla\varphi_j\dx.
 \]

 To avoid solving linear systems, the
coefficients
 $\mathbf{g}_i$ can be redefined as
convex combinations
 \[
 \mathbf{g}_i =\frac{1}{m_i}\sum_{e\in\mathcal{E}_i}m_i^e \nabla u_h|_{K_e}(\bx_i)
 \]
 of the one-sided limits $\nabla u_h|_{K_e}(\bx_i)$ in mesh cells $K_e$
 containing the point $\bx_i$. The indices of these cells are stored
 in the set $\mathcal E_i$. The weights $m_i^e\ge 0$ must add up to
 $m_i>0$. For example, the positive diagonal entries $m_i^e
 =\int_{K_e}\varphi_i\dx$ of the lumped element mass matrix can
 be used in the context of Bernstein finite element
 approximations. This definition was adopted in \cite{kuzmin2020g},
 and the corresponding operator $s_h^{e,H}(\cdot,\cdot)$
 was found to be well suited for linear
 stabilization purposes.

 The bilinear form of the VMS stabilization term
 \eqref{stab:HO} is nonsymmetric. As a consequence, its
 contribution to \eqref{semi-stab}
 may produce entropy \cite{kuzmin2020g}. The symmetric
 version
 \beq\label{stab:HOS} 
s_h^{e,H}(u_h,w_h)=\nu_e
\int_{K_e}(\nabla w_h-\mathbf{g}_h(w_h))
\cdot(\nabla u_h-\mathbf{g}_h(u_h))\dx
\eeq
is truly dissipative because it is coercive in the sense that 
$s_h^{e,H}(v_h,v_h)\ge 0\ \forall v_h\in V_h$.
Projection-based stabilization operators of this type
were proposed, for instance, in \cite{braack2006,burman2012}.
In the numerical experiments of Section~\ref{sec:num}, we use formula \eqref{stab:HOS} with
$\mathbf{g}_h$ defined by \eqref{gradrec:mcons}.

As announced in the introduction, our objective is to combine
a low-order stabilization operator and a high-order
one in order to construct a multidimensional high-order
CG-WENO version of the JST scheme \cite{jameson1981}. Introducing
a blending factor $\gamma_e\in[0,1]$, we define
(cf. \cite{barrenechea2017c})
\begin{align}
s_h^e(u_h,w_h)&=\omega\gamma_e\nu_e\int_{K_e}
(\nabla w_h-\pi_h\nabla w_h)
\cdot(\nabla u_h-\pi_h\nabla u_h)\dx,\nonumber\\
&+(1-\gamma_e)\nu_e
\int_{K_e}\nabla w_h\cdot\nabla u_h\dx.
\label{stab:WENO}
\end{align}
The additional parameter $\omega\in [0,1]$ can be used to adjust
the levels of high-order dissipation as in \cite{kuzmin2020g,lohmann2017}. 
Note that $s_h^e(u_h,w_h)$ defined by \eqref{stab:WENO} reduces to
\eqref{stab:LO} for $\gamma_e=0$ and to
\eqref{stab:HOS} for $\gamma_e=1=\omega$.

According to the JST
design philosophy, $\gamma_e$ should approach 0
in troubled cells and 1 in smooth ones. The same blending
strategy can be applied to other pairs of stabilization
operators $s_h^{e,L}$ and $s_h^{e,H}$ as long as the
latter is accuracy preserving and the former is sufficiently
dissipative. In the case $p=1$, the theoretical framework
developed by Barrenechea et al. \cite{barrenechea2017c,barrenechea2017a,barrenechea2018}
can be used to prove the validity of a discrete
maximum principle for specific choices of
$\gamma_e$ and $s_h^{e,L}$.
Numerical schemes of order $q>2$
cannot be locally bound preserving in general \cite{zhang2011}.
However, preservation of global bounds does not impose an 
order barrier.  Zhang and Shu
\cite{zhang2010b,zhang2010c,zhang2011,zhang2012} introduced
a simple scaling limiter that makes a high-order DG-WENO scheme
positivity preserving. In the
CG setting, the element-based limiters developed
in \cite{dobrev2018,kuzmin2020c,kuzmin2023,lohmann2017} can
be used to enforce positivity preservation in a similar way.

\section{Smoothness sensor}

Many shock-capturing methods and selective limiting
techniques for finite element schemes are equipped with
smoothness indicators for detection of troubled cells (see, e.g.,
\cite{hajduk2020b,krivodonova2004,lohmann2017,persson2006}).
In principle, any of these shock detectors can be used to define
$\gamma_e$ for \eqref{stab:WENO}. However, not all of the
resulting hybrid schemes will meet the conflicting demands 
for high-order accuracy and strong stability.
Motivated by the tremendous
success of WENO schemes in the context of finite volume and
DG methods, we define our dissipation-based 
stabilization term \eqref{stab:WENO} using the smoothness sensor
\begin{align}
  \gamma_e = 1 - \min\bigg(1,\frac{\|u_h^e-u_h^{e,*}\|_e}{\|u_h^e\|_e}\bigg)^q,
  \label{WENO:alpha}
\end{align}
where $u_h^e=u_h|_{K_e}$ and $u_h^{e,*}$ is 
a WENO reconstruction.
The choice of $q\ge 1$ determines how sensitive
$\gamma_e$ is to the relative difference between $u_h^e$ and $u_h^{e,*}$.
  The norm $\|\cdot\|_e$ is defined similarly
  to smoothness indicators for WENO schemes.
In this work, we use the scaled Sobolev semi-norm
  (cf. \cite{friedrich1998,jiang1996}) 
\begin{align}
  \|v\|_{e}=\left(\sum_{1\leq|\mathbf{k}|\leq p}
  h_e^{2|\mathbf{k}|-d}\int_{K_e}|D^\mathbf{k}v|^2\dx
    \right)^{1/2}\qquad
  \forall v\in H^p(K_e).
  \label{WENO:norm}
\end{align}
  In this formula, 
  $\mathbf{k}=(k_1,\ldots,k_d)$ is the
  multiindex of the partial derivative
  \[
  D^\mathbf{k}v=\pd{^{|\mathbf{k}|}v}{ x_1^{k_1}\cdots
    \partial x_d^{k_d}},\qquad |\mathbf{k}|=k_1+\ldots+k_d.
  \]

  \begin{remark}
    Derivative-based metrics of the form $\|\cdot\|_{e}^q,\ q\in\{1,2\}$ were
    used to measure the smoothness of candidate polynomials by
    Jiang and Shu \cite{jiang1996} and Friedrich \cite{friedrich1998}.
    Our definition of the blending parameter
    $\gamma_e\in[0,1]$ supports the possibility of
    raising  $\|\cdot\|_{e}$ to a user-defined power $q\ge 1$.
      \end{remark}
      
    \begin{remark}
  In contrast to WENO-based slope limiting approaches,
  we do not overwrite $u_h^e$ by $u_h^{e,*}$ in troubled cells.
  The reconstruction
  is used only to calculate $\gamma_e$ for the dissipative
  stabilization term~\eqref{stab:WENO}. The evolution of
  $u_h$ is governed by \eqref{semi-stab}. Thus the semi-discrete
  problem has the structure of a GG method, while overwriting
  DG approaches have more in common with finite volume methods.
    \end{remark}
    
\section{WENO reconstruction}

The computation of $u_h^{e,*}$ is based on the same methodology as WENO averaging of candidate polynomials in finite volume \cite{friedrich1998,jiang1996} and Runge--Kutta discontinuous Galerkin \cite{luo2007,zhong2013,zhu2009,zhu2017} methods. The main idea is to blend several approximations using a smoothness sensor to assign larger weights to less oscillatory polynomials. In the finite element context, convex combinations of Lagrange and/or Hermite polynomials corresponding to different stencils can be constructed in this way.

To facilitate the reproducibility of the numerical results to be presented in Section \ref{sec:num}, we outline the employed reconstruction procedure without claiming originality or superiority to existing alternatives. The main added value of our work is not the way to calculate the WENO polynomials $u_h^{e,*}$ but the manner in which we use them to stabilize the underlying CG scheme (see above).

Let $K_e\in\mathcal T_h$ be a generic mesh cell and $u_h$ a finite element approximation from the previous iteration, (pseudo-)time step, or Runge--Kutta stage. To measure the smoothness of $u_h^e$ on  $K_e$, we need to construct a polynomial $u_h^{e,*}$ that has the same degree $p\ge 1$ and is free of spurious oscillations. We denote by $\mathcal S^e$ the integer set containing the indices $e'$ of all  neighbor cells $K_{e'}\in\mathcal T_h$
that provide data for the computation of $u_h^{e,*}$. Note that $e\in\mathcal S^e$ by definition.
For $l=1,\ldots,m_e$, we define a stencil $\mathcal S^e_l$ as a subset of $\mathcal S^e$ and reconstruct a candidate polynomial $u_{h,l}^{e}$ from the restriction of $u_h$ to the patch $\bar\Omega_l^e=\bigcup_{e'\in\mathcal S_l^e}
K_{e'}$. The identity mapping $u_{h,0}^{e}=u_h^e$ corresponds to the stencil $\mathcal S_0^e=\{e\}$.

The properties of a reconstructed polynomial $u_{h,l}^{e}$ depend on the choice of the stencil. To achieve high-order accuracy and stability in regions where the solution is sufficiently smooth, stencils that are centered w.r.t.~$K_e$ should be included. On the other hand, one-sided stencils are needed to avoid strong oscillations around discontinuities. The number of stencils should be as small as possible to minimize the computational cost. However, it must be sufficiently large to ensure the existence of nonoscillatory candidate polynomials. An important advantage of high-order finite element methods (compared to Godunov-type finite volume schemes) is that not only cell averages but also partial derivatives of degree up to $p$ are available. Hence, direct neighbors of $K_e$ typically provide enough data.

When it comes to the construction of candidate polynomials $u_{h,l}^{e}$,
we have a choice between least-squares
fitting and interpolation. The latter approach is more common in the finite
element setting.
The cell averages\,/\,pointwise values
of Lagrange interpolation polynomials and
(some) partial derivatives of Hermite interpolation polynomials
match those of $u_h$ in cells belonging to reconstruction stencils.
The HWENO limiter developed by Luo et al. \cite{luo2007} for
DG-$\mathbb{P}_1$ approximations on simplex meshes combines $d+1$ linear Lagrange polynomials and $d+2$ linear Hermite polynomials. The reconstruction
procedure is simple and efficient, especially if the Taylor
basis \cite{luo2008} is used to represent $u_h$.
The approach
proposed by Zhong and Shu \cite{zhong2013} extends the DG polynomials
$u_h^{e'}$ of neighbor cells $e'\in\mathcal S^e_l$ into $K_e$ and
corrects the average. We use this kind of Hermite interpolation
and give some details below.

The vertex neighborhood of $K_e$ is the set of all cells that
have a common vertex with $K_e$. In our implementation, we use
only cells belonging to the von Neumann neighborhood. That is,
an index $e'$ belongs to $\mathcal S^e$ if 
$K_e$ and  $K_{e'}$ have a common boundary
(a point in 1D, an edge in 2D, a face in 3D).
Since $e\in\mathcal S^e$, we can use $\mathcal S_0^e=\{e\}$ and
$m_e\ge d+1$ reconstruction stencils $S_l^e=\{e,e'\}$,
where $e'\in\mathcal S^e\backslash\{e\}$. Using
the local basis functions of
element $e'$ and the
corresponding degrees
of freedom, we define a polynomial $u_h^{e'}$, which can
 be evaluated at $\mathbf{x}\in K_e$ to construct
(cf.~\cite{zhong2013})
\begin{align}
  u_{h,l}^e(\mathbf{x})=u_h^{e'}(\mathbf{x})+\pi_e(u_h^e-u_h^{e'}),
  \qquad\mathbf{x}\in K_e,
\label{WENO:poly}
\end{align}
where $\pi_ev=\frac{1}{|K_e|}\int_{K_e}v\dx$ denotes the 
average value of $v\in\{u_h^e,u_h^{e'}\}$ in $K_e$. The
so-defined candidate polynomials $u_{h,l}^e$ correspond
to a Hermite WENO reconstruction such that 
\[
\pi_eu_{h,l}^e= \pi_eu_h^e, \qquad
D^{\mathbf{k}}u_{h,l}^e=D^{\mathbf{k}}u_{h}^{e'},\qquad
1\le|\mathbf{k}|\le p.
\]

If the approximation
$u_h$ varies smoothly on all cells with indices in $\mathcal S^e$, then
a convex combination
\[
\tilde u_h^e=\sum_{l=0}^{m_e}\tilde{\omega}_l^eu_{h,l}^e\in\mathbb{P}_p(K_e)
\]
of candidate polynomials can be defined using \emph{linear weights}
$\tilde{\omega}_l^e \in [0,1]$. In finite volume methods, the choice
of $\tilde{\omega}_l^e$ is optimal if it maximizes the accuracy 
of reconstructed interface values at flux evaluation points. For
example, the linear component of the WENO scheme proposed by
Jiang and Shu \cite{jiang1996} yields a fifth-order approximation
that represents a convex average of three third-order approximations.
In the finite element context, small positive weights $\tilde{\omega}_l^e$
may be assigned to $u_{h,l}^e,\ l=1,\ldots,m_e$ and a large weight
$\tilde{\omega}_0^e=1-\sum_{l=1}^{m_e}\tilde{\omega}_l^e$ to
$u_{h,0}^e=u_h^e$. Such definitions can be found, e.g.,
in  \cite{zhong2013,zhu2017}. 
Since our smoothness indicator measures deviations from the WENO
reconstruction, assigning a smaller linear weight to $u_{h,0}^e$ 
results in more reliable shock detection and
stronger nonlinear stabilization.

Since $\tilde u_h^e$ may be oscillatory, the
\emph{nonlinear} weights 
${\omega}_l^e\in[0,1]$ of an adaptive WENO reconstruction 
\[
\tilde u_h^{e,*}=\sum_{l=0}^{m_e}{\omega}_l^eu_{h,l}^e\in\mathbb{P}_p(K_e)
\]
are commonly defined using linear weights $\tilde{\omega}_l^e$
and smoothness sensors $\beta_l^e$ such as
\[
\beta_l^e=\|u_{h,l}^e\|_e^q,\qquad q\ge 1,
\]
where $\|\cdot\|_e$ is the semi-norm defined by \eqref{WENO:norm}.
The classical smoothness indicators proposed by Jiang and Shu
\cite{jiang1996} and Friedrich \cite{friedrich1998} use
$q=2$ in 1D and $q=1$ in 2D, respectively.

A popular definition of the nonlinear weights $\omega_l^e$ for a WENO scheme
is given by \cite{jiang1996,zhu2009,zhu2017}
\[
\omega_l^e=\frac{\tilde{w}_l^e}{\sum_{k=0}^{m_e}\tilde{w}_k^e}, \qquad \tilde{w}_l^e = \frac{\tilde{\omega}_l^e}{(\epsilon+\beta_l^e)^r}.
\]
Here $r$ is a positive integer and $\epsilon$ is a small positive real number, which is added to avoid division by zero. In our implementation, we use the parameter settings $r=2$ and  $\epsilon = 10^{-6}$.

\begin{remark}
By definition \eqref{WENO:norm} of the semi-norm $\|\cdot\|_e$,
our WENO-based smoothness indicator \eqref{WENO:alpha} depends
only on the derivatives of candidate polynomials. The
correction of cell averages by adding $\pi_e(u_h^e-u_h^{e'})$ to
$u_h^{e'}(\mathbf{x})$ in \eqref{WENO:poly} has no influence
on the value of $\gamma_e$ and is, therefore, unnecessary in practice
(in contrast to DG-WENO schemes that overwrite $u_h^e$ by
$u_h^{e,*}$).
\end{remark}

\begin{remark}
Additional Hermite polynomials can be generated by including further stencils
of the form $\{l,m\}$ such that $K_{l}$ and $K_{m}$ share a vertex
with $K_e$. Moreover, Lagrange interpolation polynomials may be constructed.
In the context of HWENO slope limiting for piecewise-linear DG approximations,
Luo et al. \cite{luo2007} interpolate
the values of $u_h$ at the barycenters of $d+1$ cells belonging to the
reconstruction stencil $S_l^e$. In general, each mesh element has $N_K$
vertices. Let $S_l^e$ contain the indices of $N_K$ cells that surround
a vertex of $K_e$. Define $K_{e,l}$ as the convex hull of the barycenters
of these cells. Then a candidate polynomial of degree $p\ge 1$ can be
constructed by interpolating the values of $u_h$ at the nodal points of
a $p$th degree Lagrange finite element approximation on $K_{e,l}$.
This reconstruction strategy generalizes the methodology proposed
in \cite{luo2007} to arbitrary-order finite elements.
\end{remark}

In DG schemes equipped with Hermite WENO limiters, troubled cell indicators
are commonly employed to minimize the cost incurred by polynomial reconstructions
(see, e.g., \cite{zhong2013,zhu2016}). Clearly, we also have the option of skipping
the calculation of $u_{h}^{e,*}$ and setting $\gamma_e\equiv 1$ in cells
identified as smooth by an inexpensive (but safe) shock detector. To better understand
the behavior of dissipation-based WENO stabilization, we do not attempt to
localize it in this way in the present paper.

\section{Error estimates}

Following the theoretical investigations of algebraic flux correction schemes in \cite{barrenechea2016,barrenechea2018} and \cite[Chap.~4]{lohmann2019}, we perform preliminary error analysis of nonlinear dissipation-based WENO stabilization for the continuous Galerkin discretization of the linear Dirichlet boundary-value problem
\begin{subequations}\label{ibvpL}
\begin{alignat}{3}
  \mathbf{v}\cdot\nabla u+cu&=f &&\quad\mbox{in}\ \Omega,
    \label{ibvpL-pde}\\
 u&=0 &&\quad\mbox{on}\ \Gamma_-,\label{ch8:ibvp-bc}
\end{alignat}
\end{subequations}
where $\Gamma_-$ is the inflow boundary of $\Omega\subset\R^d$ and
$\mathbf{v}\in [\mathrm{Lip}(\Omega)]^d$ is a  Lipschitz-continuous
velocity field. The given functions $c\in L^\infty(\Omega)$ and
$f\in L^2(\Omega)$ determine the intensity of reactive terms.
To ensure the well-posedness of continuous and discrete problems, we make the usual assumption that 
\beq\label{divcond}
\exists c_0>0:\quad c-\frac12\nabla\cdot\mathbf{v}\ge c_0\qquad \mbox{a.e. in}\ \Omega.
\eeq
A stabilized
Galerkin discretization of \eqref{ibvpL} is given by 
\beq\label{semi-stabL}
a(u_h,w_h)+s_h(u_h; u_h,w_h)=b(w_h)\qquad\forall w_h\in V_h,
\eeq
where \cite[Chap.~2]{lohmann2019}
\begin{gather*}
a(u,w):=\int_\Omega(w\vec v\cdot\nabla u
  + cwu)\dx
  -\int_{\Gamma}wu\mins{0,\vec v\cdot\vec n}\ds,\\
  b(w):=\int_{\Omega}wf\dx.
\end{gather*}
In the integral over $\Gamma$, we denote by $\mathbf{n}$ the unit outward normal. Using $\eqref{stab:WENO}$ with $\omega=1$, we define
\begin{align*}
s_h(u_h;v_h,w_h)&=\sum_{e=1}^{E_h}\gamma_e(u_h)\nu_e
\int_{K_e}
(\nabla w_h-\pi_h\nabla w_h)
  \cdot(\nabla v_h-\pi_h\nabla v_h)\dx\\
  &+\sum_{e=1}^{E_h} (1-\gamma_e(u_h))\nu_e
\int_{K_e}\nabla w_h\cdot\nabla u_h\dx.
\end{align*}
 For any fixed
 $u_h\in V_h$, the mapping $s_h(u_h;\cdot,\cdot)$ is a symmetric bilinear
 form which is positive semi--definite and, therefore,
 satisfies the Cauchy--Schwarz
 inequality (cf.~\cite[eq.~(41)]{barrenechea2016})
 \beq\label{ch8:dh-cse}
 s_h(u_h;v_h,w_h)\le \sqrt{s_h(u_h;v_h,v_h)}
 \sqrt{s_h(u_h;w_h,w_h)}\qquad\forall v_h,w_h\in V_h.
 \eeq
 To measure the difference $e=u-u_h$ between a solution
 $u_h\in V_h$ of the discrete problem \eqref{semi-stabL}
 and an exact weak solution $u\in H^{p+1}(\Omega)$ of the continuous problem
 \eqref{ibvpL}, we will use the norms 
 \[\|e\|_a^2=c_0 \|e\|_{L^2(\Omega)}^2
+ \tfrac 12\int_{\Gamma}e^2|\vec v\cdot\vec n|\ds,\qquad
\|e\|_h^2=\|e\|_a^2+ s_h(u_h;I_he,I_he),
\]
where $I_h:C(\bar\Omega)\to V_h$ is the interpolation operator. Since
$s_h(u_h;\cdot,\cdot)$ has the same properties as the AFC stabilization operator $d_h(u_h;\cdot,\cdot)$ considered in \cite{barrenechea2016,barrenechea2018,lohmann2019}, we have 
\[
\| u - u_h\|_h  \le
\inf_{v\in V_h}\Big(
  \|u - v_h\|_h
  + \sup_{w_h \in V_h} \frac{a (u-v_h,w_h)}{\|w_h\|_h}
  + \sqrt{s_h(u_h;v_h,v_h)}\Big).  
  \]
  In essence, this result is an adaptation of Strang's first lemma to stabilized
  Galerkin schemes.
  Choosing $v_h=I_hu$ and following Lohmann \cite[Lemma 4.71]{lohmann2019}, we deduce that
\beq\label{strangest}
\| u - u_h\|_h  \le
  \|u - I_hu\|_a
  + \sup_{w_h \in V_h} \frac{a (u-I_hu,w_h)}{\|w_h\|_h}
  + \sqrt{s_h(u_h;I_hu,I_hu)}.
 \eeq
  An $\mathcal O(h^p)$ bound for the sum of the
    first two terms on the right-hand of this inequality can be
    derived as in \cite{hartmann2008,lohmann2019}.
    The last term measures the consistency
    error due to nonlinear stabilization.

    Recall
    that $\nu_e=\mathcal O(h)$ and $\gamma_e\in[0,1]$.
  If the exact weak solution $u$ is sufficiently smooth, then
  \[\sum_{e=1}^{E_h}\gamma_e(u_h)\nu_e\int_{K_e}
  |\nabla I_hu-\pi_h\nabla I_hu|^2\dx\le Ch \sum_{e=1}^{E_h}\int_{K_e}|\nabla I_hu|^2\dx
 =  Ch\|\nabla I_hu\|_{L^2(\Omega)}^2.
 \]
 The low-order component is bounded similarly.
 Hence,  $\sqrt{s_h(u_h;I_hu,I_hu)}=\mathcal O(h^{1/2})$ in the worst case.

 On the other hand, suppose that the high-order component
 dominates and the estimate
 \begin{align*}
  \sum_{e=1}^{E_h}(1-\gamma_e(u_h))\nu_e
  \int_{K_e}|\nabla I_hu|^2\dx&\le Ch\sum_{e=1}^{E_h}
  \int_{K_e}|\nabla I_hu-\pi_h\nabla I_hu|^2\dx
 \\ &=Ch\|\nabla I_hu-\pi_h\nabla I_hu\|_{L^2(\Omega)}^2
  \end{align*}
  holds \emph{a posteriori}
  for a specific choice of $\gamma_e$. Using the triangle inequality, we find that
  \[
  \|\nabla I_hu-\pi_h\nabla I_hu\|_{L^2(\Omega)}\le
  \|\nabla I_hu-\nabla u\|_{L^2(\Omega)}+ \|\nabla u-\pi_h\nabla u\|_{L^2(\Omega)}
  +\|\pi_h(\nabla u-I_h\nabla u)\|_{L^2(\Omega)},\quad
  \]
  where
  \[
  \|\nabla I_hu-\nabla u\|_{L^2(\Omega)}=|I_hu-u|_{H^1(\Omega)}
  \le Ch^{p}|u|_{H^{p+1}(\Omega)}.
  \]
  If $\pi_h$ is the $L^2$ projection operator, then
  its best approximation property implies that
  \[
  \|\pi_h\nabla u-
  \nabla u\|_{L^2(\Omega)}\le  \|I_h\nabla u-\nabla u\|_{L^2(\Omega)}\le
Ch^p|\nabla u|_{H^{p}(\Omega)}.
 \]
 Since $\|\pi_hv\|_{L^2(\Omega)}^2=
 (v,\pi_hv)_{L^2(\Omega)}$ for $v\in L^2(\Omega)$,
 we can use the Cauchy--Schwarz inequality
 and Young's inequality to show that $\|\pi_hv\|_{L^2(\Omega)}\le 
 \|v\|_{L^2(\Omega)}$ for
 any $v\in L^2(\Omega)$. Combining the above auxiliary results, we conclude  that
 $\sqrt{s_h(u_h;I_hu,I_hu)}=\mathcal O(h^{p+1/2})$ under our assumption that the high-order
 stabilization is dominant for the particular choice of $\gamma_e(u_h)$.
 
\begin{remark}
  An optimal $\mathcal O(h^{p+1/2})$ \emph{a priori} error estimate for the linear high-order
  scheme  can be obtained following the analysis
  of the SUPG method in
  \cite{hartmann2008,kuzmin2023}. Instead of $\tau h\|\mathbf{v}\cdot\nabla
 \rho\|_{L^2(\Omega)}$, where $\rho=u-I_hu$ is the interpolation error and
  $\tau$ is the SUPG stabilization
  parameter, the estimate will contain $\sqrt{s_h(u_h;I_hu,I_hu)}$ with
  $\gamma_e\equiv 1$. This term is $\mathcal O(h^{p+1/2})$, as shown above.
  
\end{remark}
 
 \begin{remark}
 If the restriction of $u_h$ to the union of mesh cells with indices in
 $\mathcal S^e$ is a polynomial of degree up to $p$, then
$\|u_h^e-u_h^{e,*}\|_e=0$ and $\gamma_e(u_h)=1$
 for any choice of the
 nonlinear weights. Therefore, polynomial exact
 solutions can be reproduced exactly by our method.
\end{remark}

 For a general smooth function $u\in C^{p+1}(\bar\Omega)$, we claim that
 the WENO consistency error satisfies
 \beq\label{claim}
 \sqrt{s_h(I_hu;I_hu,I_hu)}=\mathcal O(h^{p+1/2}).
 \eeq
 To prove the validity of
 this claim, we need to estimate $\|v_h^e-v_h^{e,*}\|_e$ for $v_h=I_hu$.
 The multivariate Taylor expansions of
 the Hermite candidate polynomials $v^e_{h,l}$ are given by
 $$
v^e_{h,l}(\mathbf{x})=v^e_l+\sum_{1\leq|\mathbf{m}|\leq p}
D^\mathbf{m}u(\mathbf{x}^e_{l})\frac{(\mathbf{x}-\mathbf{x}^e_{l})^\mathbf{m}}{\mathbf{m}!},
 $$
where we use the multiindex notation. For $e'\in\mathcal S_l^e$, we denote by
$\mathbf{x}^e_{l}$ the centroid $\pi_{e'}\mathbf{x}$ of $K_{e'}$. The
value of the constant $v^e_l=v^e_{h,l}(\mathbf{x}^e_{l})$
is uniquely determined by the requirement
that $v^e_{h,l}$ have the cell average $\pi_ev^e_{h,l}=\pi_ev^e_{h}=\pi_eI_hu$.
By definition of  $\|\cdot\|_e$ and $v_h^{e,*}$, we have
\beq\label{enormest}
\|v_h^e-v_h^{e,*}\|_e^2
=\sum_{1\leq|\mathbf{k}|\leq p}h_e^{2|\mathbf{k}|-d}
\int_{K_e}\left|\sum_{l=1}^{m_e}\omega_l^eD^\mathbf{k}(v^e_{h,l}-v^e_{h,0})\right|^2\dx.
\eeq
Introducing the Taylor basis functions \cite{luo2008}
$$
\psi^e_{\mathbf{m}}(\mathbf{x})=\frac{1}{\mathbf{m}!}\left[
(\mathbf{x}-\mathbf{x}^e_{0})^{\mathbf{m}}-\frac{1}{|K_e|}\int_{K_e}
  (\mathbf{x}-\mathbf{x}^e_{0})^{\mathbf{m}}\dx\right]
$$
such that
$$
v^e_{h,l}(\mathbf{x})=\pi_ev^e_{h,l}+\sum_{1\leq|\mathbf{m}|\leq p}
 D^\mathbf{m}v^e_{h,l}(\mathbf{x}^e_{0})\psi^e_{\mathbf{m}}(\mathbf{x}),
$$
we proceed to estimate 
$$
D^\mathbf{k}v^e_{h,l}(\mathbf{x})-D^\mathbf{k}v^e_{h,0}(\mathbf{x})
=\left(
\sum_{l=1}^{m_e}\omega_l^e
\sum_{1\leq|\mathbf{m}|\leq p}
    [D^\mathbf{m}v^e_{h,l}(\mathbf{x}^e_{0})-
      D^\mathbf{m}u(\mathbf{x}^e_{0})]\right)D^\mathbf{k}
\psi^e_{\mathbf{m}}(\mathbf{x}).
$$

By definition of the candidate polynomials $v^e_{h,l},\ l\in\{1,\ldots,m_e\}$, the difference
$$
R^e_{\mathbf{m},l}(\mathbf{x})=D^\mathbf{m}v^e_{h,l}(\mathbf{x})-D^\mathbf{m}u(\mathbf{x})
$$
is the remainder of a truncated Taylor expansion of $D^\mathbf{m}u(\mathbf{x})$
about the centroid $\mathbf{x}^e_{l}$. The degree of the Taylor
polynomial $D^\mathbf{m}v^e_{h,l}(\mathbf{x}^e_{0})$ is $p-|\mathbf{m}|\ge 0$.
It follows that $$|R^e_{\mathbf{m},l}(\mathbf{x}^e_{0})|=\mathcal O(h^{p-|\mathbf{m}|+1}).$$
Invoking the definition of $\psi^e_{\mathbf{m}}$, we deduce that
$|D^\mathbf{k}
\psi^e_{\mathbf{m}}(\mathbf{x})|=\mathcal O(h^{|\mathbf{m}|-|\mathbf{k}|})$
for $|\mathbf{x}-\mathbf{x}^e_{0}|\le Ch$. Thus
$$|
D^\mathbf{k}v^e_{h,l}(\mathbf{x})-D^\mathbf{k}v^e_{h,0}(\mathbf{x})|
=\mathcal O(h^{p-|\mathbf{k}|+1}).
$$
Using this result and the fact that $|K_e|=\mathcal O(h^d)$
to estimate the right-hand side of 
\eqref{enormest}, we arrive at
$$
\|v_h^e-v_h^{e,*}\|_e^2=\mathcal O(h^{2p+2}).
$$

Let the smoothness indicator $\gamma_e(v_h)$ be defined by formula \eqref{WENO:alpha} with $q=2$. Then
\begin{align*}
  (1-\gamma_e(v_h))\nu_e\int_{K_e}|v_h|^2\dx&=
  \nu_e\mina{1,\frac{\|v_h^e-v_h^{e,*}\|_e^2}{\|v_h^e\|_e^2}}\|v_h^e\|_{L^2(\Omega)}^2
\\ &\le C\nu_eh^{d-2}\mina{\|v_h^e-v_h^{e,*}\|_e^2,\|v_h^e\|_e^2},
\end{align*}
where $\nu_e=\mathcal O(h)$ and $\|v_h^e\|_e=\mathcal O(1)$. Hence,
the low-order component of $s_h(I_hu;I_hu,I_hu)$ is
$\mathcal O(h^{2p+1})$. The corresponding estimate for the high-order
component was obtained above. This proves the validity of our claim
that the consistency error of the WENO stabilization satisfies \eqref{claim}.

The nonlinear term that appears in \eqref{strangest} is not
$s_h(I_hu;I_hu,I_hu)$ but $s_h(u_h;I_hu,I_hu)$. To obtain a final
\emph{a priori} error estimate, let us assume that the Lipschitz
continuity condition
$$
|s_h(u_h;v_h,v_h)-s_h(v_h;v_h,v_h)|\le s_h(u_h;v_h-u_h,v_h-u_h)\qquad \forall u_h\in V_h
$$
holds for the interpolant $v_h=I_hu$ of $u\in C^{p+1}(\bar\Omega)$. Then \eqref{strangest}
implies
\begin{align*}
  \| u - u_h\|_a &=\| u - u_h\|_h-\sqrt{s_h(u_h;I_hu-u_h,I_hu-u_h)}
   \\&\le
  \|u - I_hu\|_a
  + \sup_{w_h \in V_h} \frac{a (u-I_hu,w_h)}{\|w_h\|_h}
  -\sqrt{s_h(u_h;I_hu-u_h,I_hu-u_h)}\\
  &+ \sqrt{s_h(u_h;I_hu-u_h,I_hu-u_h)+s_h(I_hu;I_hu,I_hu)}
  \\&\le
  \|u - I_hu\|_a
  + \sup_{w_h \in V_h} \frac{a (u-I_hu,w_h)}{\|w_h\|_h}
  + \sqrt{s_h(I_hu;I_hu,I_hu)}.
\end{align*}
  The  last term on the right-hand side of this inequality
  can be now be estimated using \eqref{claim}. As already mentioned,
  an $\mathcal O(h^p)$ estimate is available for 
  the sum of the first
  two terms.

  \begin{remark}
    Barrenechea et al. \cite{barrenechea2018} use similar arguments to
    derive an improved \emph{a priori} error estimate for algebraic
    flux correction schemes equipped with linearity-preserving limiters.
    \end{remark}

\section{Numerical examples}
\label{sec:num}


In this section, we present the results of some numerical experiments for linear and nonlinear scalar problems. The objective of our study is to show that the CG-WENO scheme
\begin{itemize}
\item exhibits optimal convergence behavior for problems with smooth solutions; 
\item introduces sufficient numerical dissipation in the neighborhood of shocks;
\item is unlikely to produce entropy-violating solutions in the nonlinear case.
\end{itemize}
In our description of the numerical results, we label the methods under
investigation as follows:
\medskip

\begin{tabular}{c|l}
  CG& continuous Galerkin method without any stabilization;\\
 VMS& CG + two-level variational multiscale stabilization \eqref{stab:HO};\\
 HO& CG + linear high-order stabilization, i.e., \eqref{stab:WENO}
  using $\gamma_e\equiv 1$;\\
 LO& CG + linear low-order stabilization, i.e., \eqref{stab:WENO}
  using $\gamma_e\equiv 0$;\\
 WENO& CG + nonlinear stabilization \eqref{stab:WENO}
      using $\gamma_e$ defined by \eqref{WENO:alpha}.
\end{tabular}
\medskip

The consistent-mass $L^2$ projection operator $\pi_e$ is used to calculate
 $\mathbf{g}_h=\pi_e\nabla u_h$
in~\eqref{stab:HO} and the projected gradients in~\eqref{stab:WENO}.
We set  $\omega=1.0$ in \eqref{stab:WENO} and $q=1$ in~\eqref{WENO:alpha}
unless mentioned otherwise. The default setting for the linear weights
of the Hermite WENO reconstruction is (cf. \cite{zhu2017})
$$\tilde \omega_0^e=1-m_e\cdot 10^{-3},\qquad
\tilde{\omega}_l^e=10^{-3},\quad l=1,\ldots,m_e.
$$

Computations are performed using $\mathbb{Q}_p$ Lagrange finite elements of degree $p=\{1,2,3,4\}$ and explicit Runge--Kutta methods of order $p+1$. For comparison purposes, we present numerical solutions corresponding to different values of $p$ for a fixed total number of degrees of freedom (DoFs), which we denote by $N_h$. That is, we coarsen the mesh if we increase the polynomial degree.
The time step $\Delta t$ is chosen to be proportional to the mesh size
$h$ and small enough for the temporal error to be negligible. Spatial
discretization errors are measured using the $L^1$ norm.
The experimental order of convergence (EOC)
is determined as in \cite{leveque1996,lohmann2017}. In captions of some
figures we use the shorthand notation $E_1$ for $\|u_h-u_{\text{exact}}\|_{L^1(\Omega)}$,
where $u_{\text{exact}}$ is the exact solution.
The implementation of all schemes is based on the open-source C++  finite element library MFEM \cite{anderson2021,mfem}.

\subsection{One-dimensional linear advection with constant velocity in 1D}
The first problem we consider is the one-dimensional linear advection equation
\begin{align}
\frac{\partial u}{\partial t}+v\frac{\partial u}{\partial x} = 0 \quad \mbox{in}\ \Omega=(0,1)
\label{num:adveq}
\end{align}
with constant velocity $v=1$. To begin, we advect the smooth initial condition
\begin{align}
u_0(x)=\cos(2\pi(x-0.5))
\label{num:advinitone}
\end{align}
up to the final time $t=1.0$ using $q=3$ in the formula for $\gamma_e$. The results of a grid convergence study on uniform meshes are shown in Table \ref{tab:convlin}. All methods deliver the optimal $L^1$ convergence rates (EOC$\,\approx p+1$) for Lagrange finite elements of degree $p\in\{1,2,3,4\}$.

\begin{table}[h!]
	\begin{subtable}{\textwidth}
		\centering
		\begin{tabular}{||c||c|c||c|c||c|c||}
			\hline
			&\multicolumn{2}{c||}{CG}&\multicolumn{2}{c||}{VMS}&\multicolumn{2}{c||}{WENO}\\
			\hline
			$N_h$ & $\|u_h-u_{\text{exact}}\|_{L^1}$ & EOC & $\|u_h-u_{\text{exact}}\|_{L^1}$ & EOC & $\|u_h-u_{\text{exact}}\|_{L^1}$ & EOC\\
			\hline
			16	&   8.20\mbox{e-3}&	2.04&   9.22\mbox{e-3}&	2.19&   9.26\mbox{e-2}&	1.98\\
			32	&   2.05\mbox{e-3}&	2.00&   2.17\mbox{e-3}&  2.09&   2.68\mbox{e-2}&  1.79\\
			64	&   5.11\mbox{e-4}&  2.00&   5.27\mbox{e-4}&  2.04&   3.24\mbox{e-3}&  3.05\\
			128 &   1.28\mbox{e-4}&  2.00&   1.30\mbox{e-4}&  2.02&   2.38\mbox{e-4}&  3.77\\
			256 &   3.20\mbox{e-5}&  2.00&   3.22\mbox{e-5}&  2.01&   3.51\mbox{e-5}&  2.76\\
			512 &   7.99\mbox{e-6}&  2.00&   8.03\mbox{e-6}& 	2.00&	8.31\mbox{e-6}&  2.08\\
			1024&   2.02\mbox{e-6}&  1.99&	2.02\mbox{e-6}& 	1.99&	2.06\mbox{e-6}&  2.01\\
			\hline
		\end{tabular}
		\caption{$p=1$}
	\end{subtable}
	\begin{subtable}{\textwidth}
		\centering
		\begin{tabular}{||c||c|c||c|c||c|c||}
			\hline
			&\multicolumn{2}{c||}{CG}&\multicolumn{2}{c||}{VMS}&\multicolumn{2}{c||}{WENO}\\
			\hline
			$N_h$ & $\|u_h-u_{\text{exact}}\|_{L^1}$ & EOC & $\|u_h-u_{\text{exact}}\|_{L^1}$ & EOC & $\|u_h-u_{\text{exact}}\|_{L^1}$ & EOC\\
			\hline
			32  &   8.17\mbox{e-4}&  1.45&   5.17\mbox{e-4}&    2.30&   2.72\mbox{e-4}&  3.09\\
			64  &   5.80\mbox{e-5}&  3.82&   9.14\mbox{e-5}&    2.50&   3.26\mbox{e-5}&  3.06\\
			128 &   4.25\mbox{e-6}&  3.77&   1.34\mbox{e-5}&    2.78&   4.06\mbox{e-6}&  3.01\\
			256 &   4.29\mbox{e-7}&  3.31&   1.75\mbox{e-6}&    2.93& 	5.06\mbox{e-7}&  3.00\\
			512 & 	5.36\mbox{e-8}&  3.00&	 2.22\mbox{e-7}&	2.98&	6.31\mbox{e-8}&  3.00\\
			1024&  	6.70\mbox{e-9}&  3.00&	 2.78\mbox{e-8}& 	3.00&	7.90\mbox{e-9}&  3.00\\
			\hline
		\end{tabular}
		\caption{$p=2$}
	\end{subtable}
	\begin{subtable}{\textwidth}
		\centering
		\begin{tabular}{||c||c|c||c|c||c|c||}
			\hline
			&\multicolumn{2}{c||}{CG}&\multicolumn{2}{c||}{VMS}&\multicolumn{2}{c||}{WENO}\\
			\hline
			$N_h$ & $\|u_h-u_{\text{exact}}\|_{L^1}$ & EOC & $\|u_h-u_{\text{exact}}\|_{L^1}$ & EOC & $\|u_h-u_{\text{exact}}\|_{L^1}$ & EOC\\
			\hline
			48  &   4.15\mbox{e-6}&   4.10&   3.10\mbox{e-6}&   4.13&   2.59\mbox{e-5}&  	3.40\\
			96  &   2.79\mbox{e-7}&   3.89&   1.86\mbox{e-7}&   4.05&   2.26\mbox{e-6}&  	3.52\\
			192 &   1.75\mbox{e-8}&   3.99&   1.16\mbox{e-8}&   4.00&   1.54\mbox{e-7}&     3.88\\
			384 &   1.09\mbox{e-9}&   4.01&   7.26\mbox{e-10}&  4.00&	9.97\mbox{e-9}&     3.95\\
			768 &   6.69\mbox{e-11}&  4.03&	  4.54\mbox{e-11}&	4.00&	6.36\mbox{e-10}&	3.97\\
			1536&   3.93\mbox{e-12}&  4.09&	  2.84\mbox{e-12}&	4.00&	4.05\mbox{e-11}&	3.97\\
			\hline
		\end{tabular}
		\caption{$p=3$}
	\end{subtable}
	\begin{subtable}{\textwidth}
		\centering
		\begin{tabular}{||c||c|c||c|c||c|c||}
			\hline
			&\multicolumn{2}{c||}{CG}&\multicolumn{2}{c||}{VMS}&\multicolumn{2}{c||}{WENO}\\
			\hline
			$N_h$ & $\|u_h-u_{\text{exact}}\|_{L^1}$ & EOC & $\|u_h-u_{\text{exact}}\|_{L^1}$ & EOC & $\|u_h-u_{\text{exact}}\|_{L^1}$ & EOC\\
			\hline
			64 &    4.95\mbox{e-7}&   3.93&   1.92\mbox{e-7}&   4.52&   6.03\mbox{e-6}&  5.46\\
			128&	9.06\mbox{e-9}&   5.77&   7.07\mbox{e-9}&   4.76&   9.35\mbox{e-8}&  6.01\\
			256&    1.63\mbox{e-10}&  5.80&   2.33\mbox{e-10}&  4.92&   1.46\mbox{e-9}&  6.00\\
			512&    3.47\mbox{e-12}&  5.55&   7.39\mbox{e-12}&  4.98&	4.30\mbox{e-11}&  5.08\\
			\hline
		\end{tabular}
		\caption{$p=4$}
	\end{subtable}
	\caption{1D linear advection, grid convergence history for finite elements of degree $p\in\{1,2,3,4\}$.}
	\label{tab:convlin}
\end{table}

\begin{figure}[!htb]
	\centering
	\begin{subfigure}[b]{.9\linewidth}
		\centering
		\begin{tikzpicture}
		\draw[rounded corners] (0, 0) rectangle (12, 0.5) node[pos=.5]{};
		\draw[very thick, color={rgb:red,1;green,0.2;blue,0.2}] (0.75,0.25)--(1.25,0.25);
		\node at (2,0.25) (a) {$p=1$};
		\node at (5,0.25) (a) {$p=2$};
		\node at (8,0.25) (a) {$p=4$};
		\node at (11,0.25) (a) {Exact};
		\draw[very thick,color={rgb:red,0.2;green,0.8;blue,0.8}] (3.75,0.25)--(4.25,0.25);
		\draw[very thick,color={rgb:red,0;green,0;blue,1}] (6.75,0.25)--(7.25,0.25);
		\draw[dashed, very thick] (9.75,0.25)--(10.25,0.25);
		\end{tikzpicture}
	\end{subfigure}
	\begin{subfigure}[b]{.4\linewidth}
		\includegraphics[width=\linewidth]{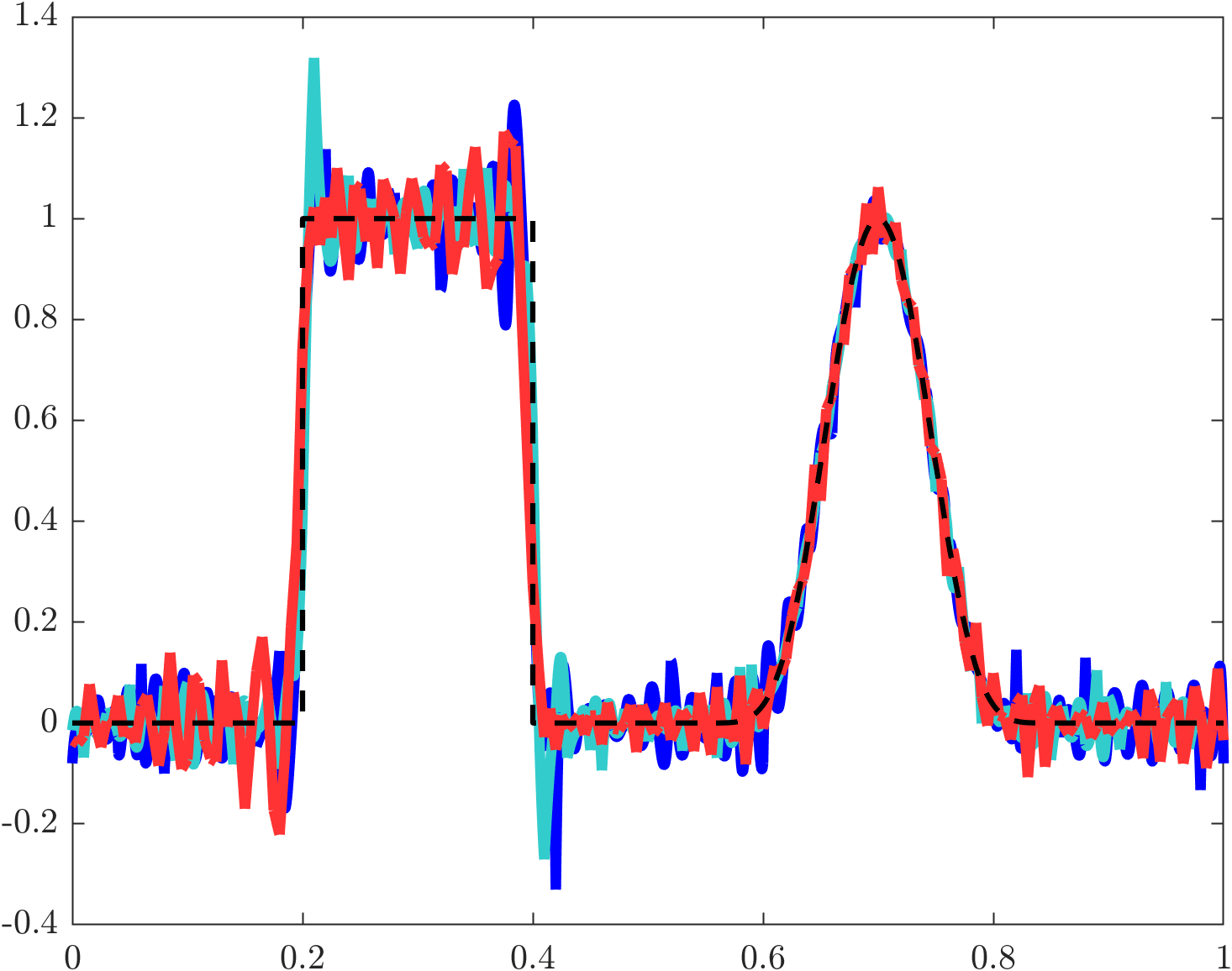}
		\caption{CG}
		\label{fig:advcg}
	\end{subfigure}
	\begin{subfigure}[b]{.4\linewidth}
		\includegraphics[width=\linewidth]{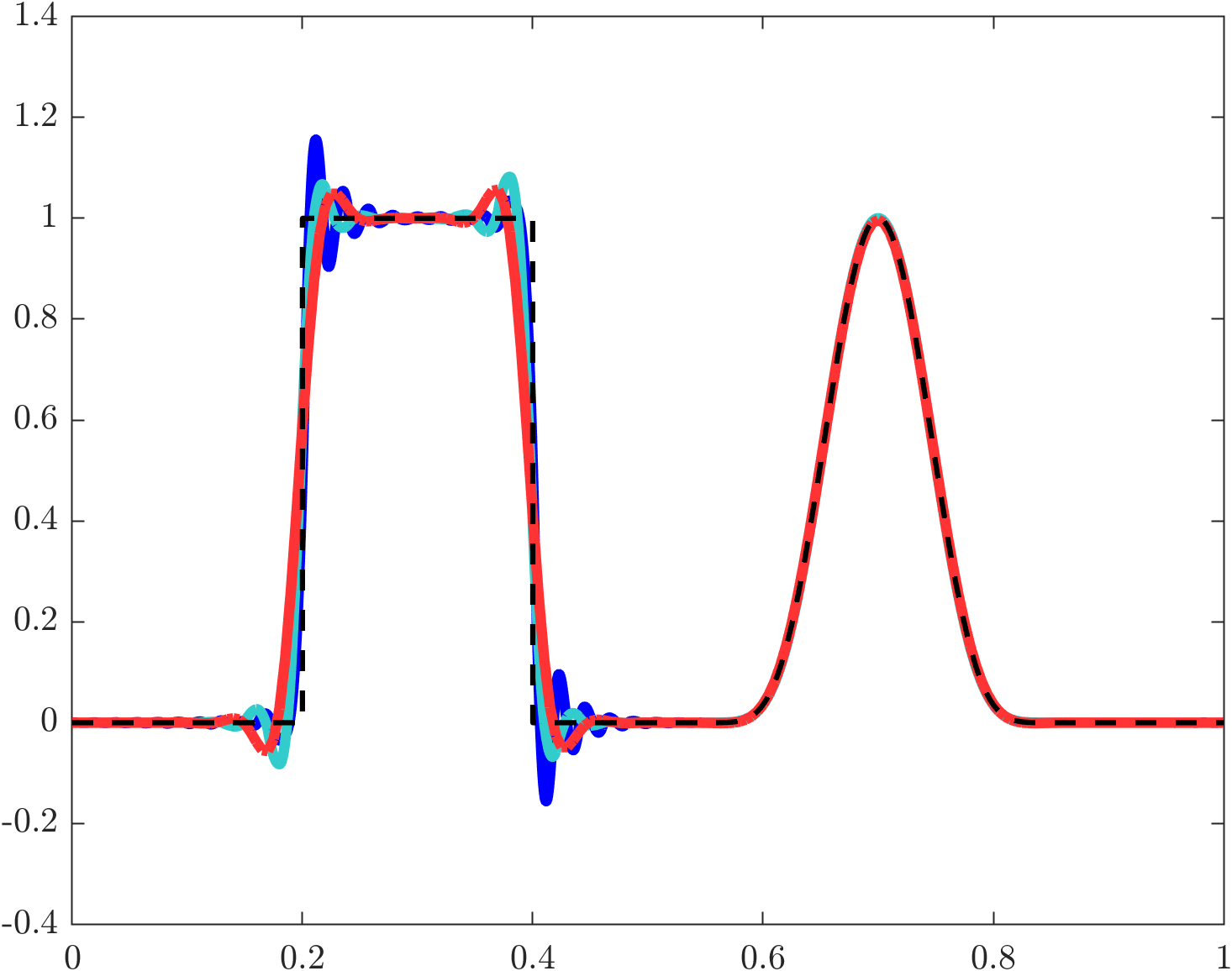}
		\caption{HO}
		\label{fig:advho}
	\end{subfigure}
	\begin{subfigure}[b]{.4\linewidth}
		\includegraphics[width=\linewidth]{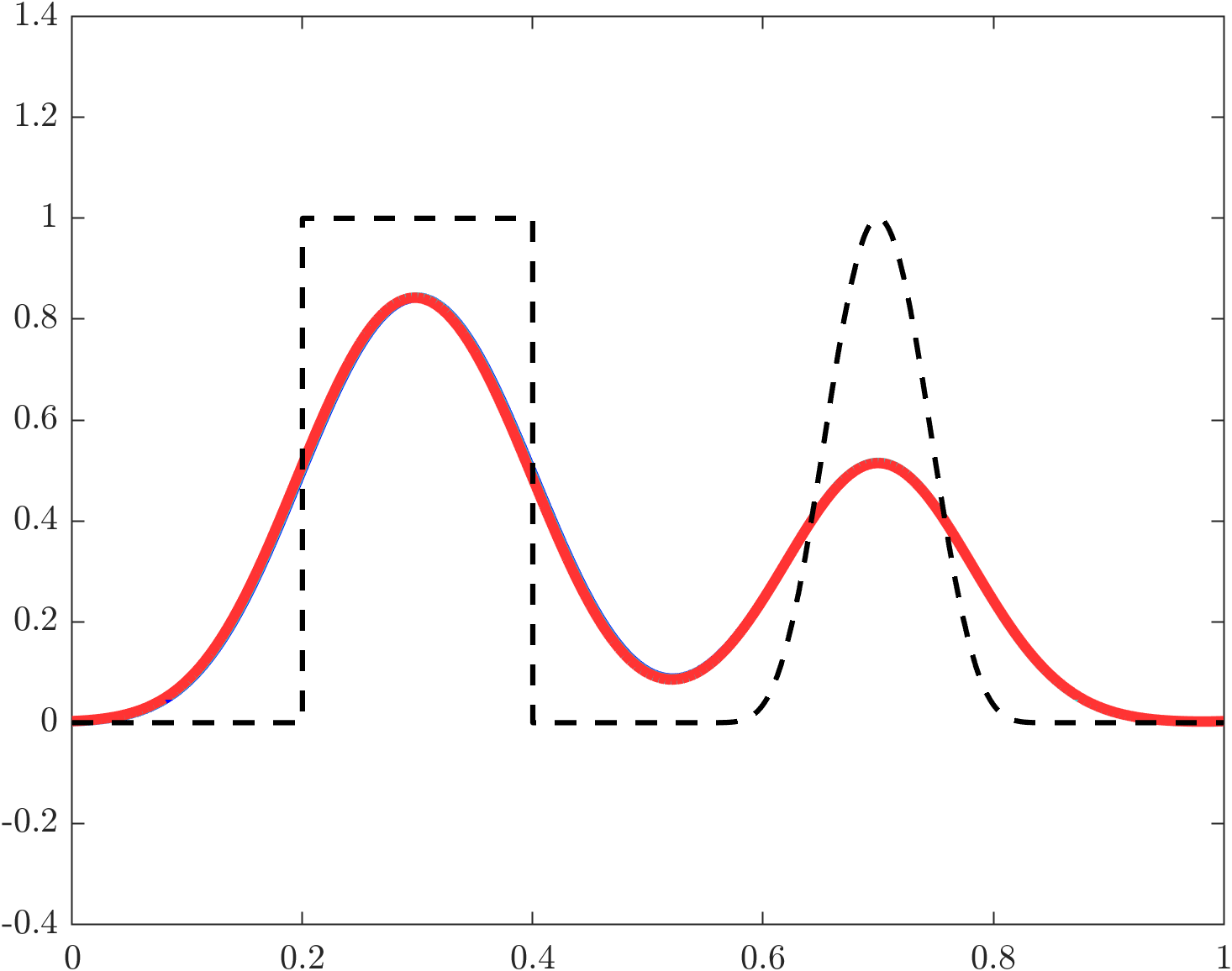}
		\caption{LO}
		\label{fig:advlo}
	\end{subfigure}
	\begin{subfigure}[b]{.4\linewidth}
		\includegraphics[width=\linewidth]{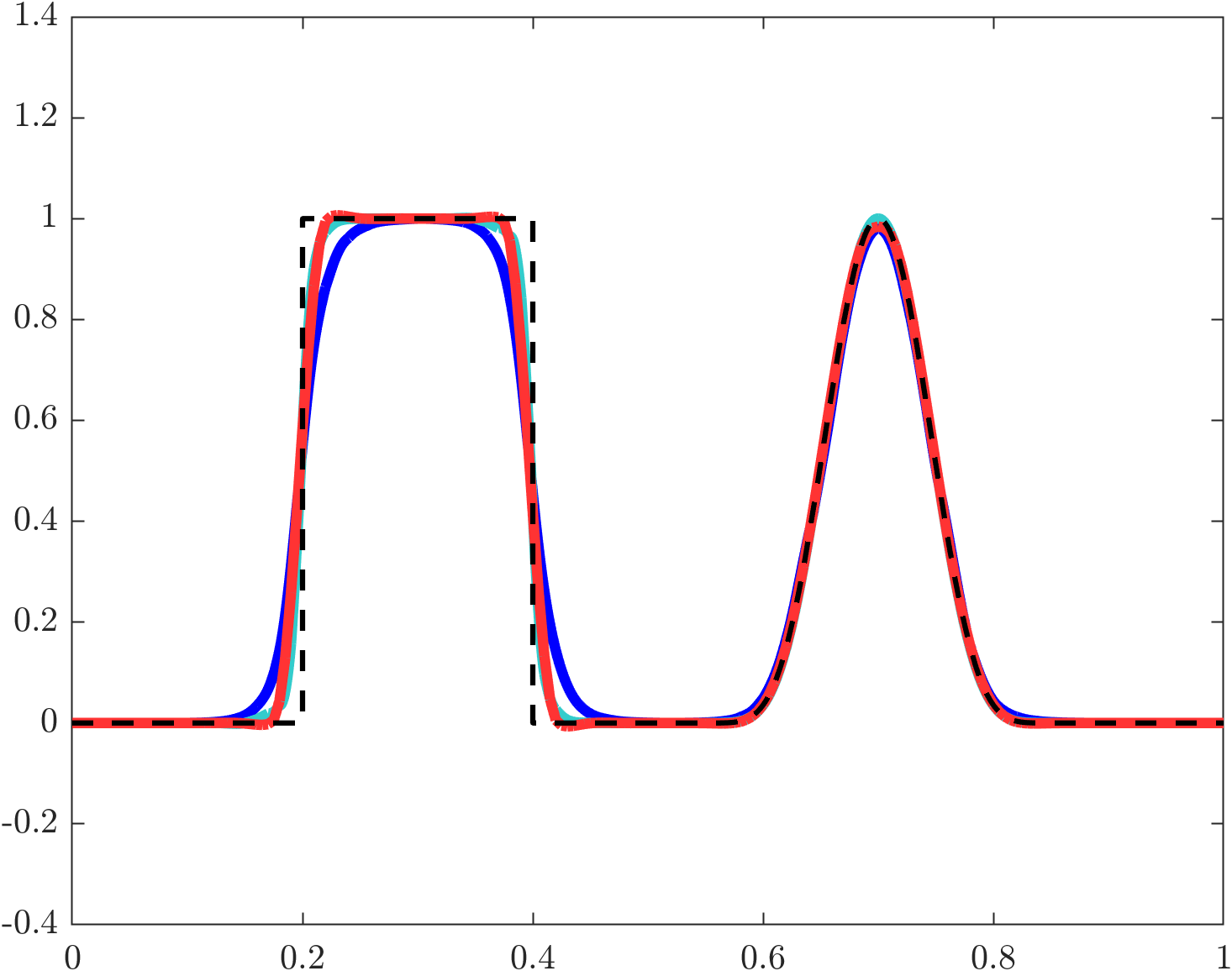}
		\caption{WENO}
		\label{fig:advweno}
	\end{subfigure}
	\caption{1D linear advection, numerical solutions at $t=1$ obtained using $N_h=200$ and $p=\{1,2,4\}$.}
	\label{fig:linadv}
\end{figure}

To study the stability properties of our schemes, we replace the above
 initial condition with \cite{hajduk2021}
\begin{align}
u_0(x) = \begin{cases}
1 & \mbox{if} \ 0.2 \le x \le 0.4, \\
\exp(10)\exp(\frac{1}{0.5-x})\exp(\frac{1}{x-0.9}) & \mbox{if}  \ 0.5 < x < 0.9, \\
0 & \mbox{otherwise}
\end{cases}
\label{num:advinittwo}
\end{align}
and choose the mesh size $h$ corresponding to $N_h=200$ DoFs for a given polynomial degree $p$.
The results for $p\in\{1,2,4\}$
are shown in Figs \ref{fig:advcg}-\ref{fig:advweno}. Additionally, we list the global
maxima and minima of numerical solutions $u_h$ in Table~\ref{table-maxmin}.
As expected, the standard CG approximation oscillates heavily. The activation of linear HO stabilization localizes spurious oscillations to a small neighborhood of discontinuities and reduces the magnitude of undershoots/overshoots. The numerical solutions produced by the linear LO scheme are completely free of oscillations. However, they are strongly dissipative because first-order artificial diffusion is added everywhere. The result shown in Fig.~\ref{fig:advweno} demonstrates the superb shock-capturing capabilities of our CG-WENO scheme. Both smooth and discontinuous portions of the advected initial profile are preserved very well compared to other methods. Stronger smearing of discontinuities for larger values of $p$ is due to the fact that we keep $N_h$ fixed.

\begin{table}[htb!]
	\centering
	\begin{tabular}{||c||c|c||c|c||c|c||c|c||}
		\hline
		&\multicolumn{2}{c||}{CG}&\multicolumn{2}{c||}{HO}&\multicolumn{2}{c||}{LO}&
		\multicolumn{2}{c||}{WENO}\\
		\cline{2-9}
		& $u_h^{\min}$ & $u_h^{\max}$  & $u_h^{\min}$ & $u_h^{\max}$& $u_h^{\min}$ & $u_h^{\max}$  & $u_h^{\min}$ & $u_h^{\max}$\\
		\hline
		$p=1$	&	-0.2210&	1.1725&	-0.0564&	1.0564&	0.0026&	0.8424&	-0.0066&	1.0066 \\
		$p=2$	&   -0.2704&	1.3186& -0.0821&	1.0821&	0.0022&	0.8428&	-0.0013&	1.0013\\
		$p=4$	&   -0.3303&	1.2257& -0.1557&	1.1560& 0.0000& 0.8427&	\phantom{-}0.0000&	0.9999\\
		\hline
	\end{tabular}
	\caption{1D linear advection, global maximum $u_h^{\max}$ and minimum  $u_h^{\min}$ of $u_h$ for $p\in\{1,2,3,4\}$.}
	\label{table-maxmin}
\end{table}

\subsection{One-dimensional inviscid Burgers equation}

In the next experiment, we solve the one-dimensional inviscid Burgers equation
\begin{align}
\frac{\partial u}{\partial t}+\frac{\partial(u^2/2)}{\partial x} = 0 \quad \mbox{in}\ \Omega=(0,1).
\label{num:bureq}
\end{align}
The initial condition for this nonlinear test problem is given by
\begin{align}
u_0(x)=\sin(2\pi x).
\label{num:burinit}
\end{align}
It can be shown that the unique entropy solution develops a shock at the critical time $t_c=\frac{1}{2\pi}$. In our grid convergence studies, we stop computations at $t=0.1 < t_c$, because this final time is small enough for the exact solution to remain sufficiently smooth. Table~\ref{tab:convbur} shows that optimal EOCs can again be achieved using HO and WENO stabilization. The $L^1$ convergence rate of the standard CG method is also as high as $p+1$ for odd polynomial degrees but drops down to $p$ for even ones.

\begin{table}[h!]
	\begin{subtable}{\textwidth}
		\centering
		\begin{tabular}{||c||c|c||c|c||c|c||}
			\hline
			&\multicolumn{2}{c||}{CG}&\multicolumn{2}{c||}{VMS}&\multicolumn{2}{c||}{WENO}\\
			\hline
			$N_h$ & $\|u_h-u_{\text{exact}}\|_{L^1}$ & EOC & $\|u_h-u_{\text{exact}}\|_{L^1}$ & EOC & $\|u_h-u_{\text{exact}}\|_{L^1}$ & EOC\\
			\hline
			16  & 	1.22\mbox{e-2}&	 2.13&   1.19\mbox{e-2}&	2.13&	2.16\mbox{e-2}&	2.24\\
			32  & 	2.56\mbox{e-3}&  2.25&   2.55\mbox{e-3}&	2.22&	5.81\mbox{e-3}&	1.90\\
			64  & 	6.05\mbox{e-4}&  2.08&   6.04\mbox{e-4}& 	2.08&	1.10\mbox{e-3}&	2.40\\
			128 & 	1.49\mbox{e-4}&  2.02&   1.49\mbox{e-4}&    2.02&	1.56\mbox{e-4}&	2.82\\
			256 & 	3.72\mbox{e-5}&  2.00& 	 3.72\mbox{e-5}& 	2.00&	3.71\mbox{e-5}&	2.07\\
			512 &  	9.29\mbox{e-6}&  2.00& 	 9.29\mbox{e-6}&	2.00&	9.28\mbox{e-6}&	2.00\\
			1024&  	2.33\mbox{e-6}&  2.00& 	 2.33\mbox{e-6}& 	2.00&	2.32\mbox{e-6}&	2.00\\
			\hline
		\end{tabular}
		\caption{$p=1$}
	\end{subtable}
	\begin{subtable}{\textwidth}
		\centering
		\begin{tabular}{||c||c|c||c|c||c|c||}
			\hline
			&\multicolumn{2}{c||}{CG}&\multicolumn{2}{c||}{VMS}&\multicolumn{2}{c||}{WENO}\\
			\hline
			$N_h$ & $\|u_h-u_{\text{exact}}\|_{L^1}$ & EOC & $\|u_h-u_{\text{exact}}\|_{L^1}$ & EOC & $\|u_h-u_{\text{exact}}\|_{L^1}$ & EOC\\
			\hline
			32  &  	8.51\mbox{e-4}&  2.99& 	8.07\mbox{e-4}&  3.03&	5.70\mbox{e-4}&	3.27\\
			64  &  	2.73\mbox{e-4}&  1.64& 	2.11\mbox{e-4}&	1.94&	1.47\mbox{e-4}&	1.95\\
			128 &  	6.34\mbox{e-5}&  2.11&	3.44\mbox{e-5}&	2.62&	1.91\mbox{e-5}&	2.94\\
			256 & 	1.52\mbox{e-5}&  2.06&	5.21\mbox{e-6}&	2.72&	2.33\mbox{e-6}&	3.04\\
			512 &  	3.73\mbox{e-6}&  2.02&	7.68\mbox{e-7}&	2.76&	2.87\mbox{e-7}&	3.02\\
			1024& 	9.27\mbox{e-7}&  2.01&	1.10\mbox{e-7}&	2.80&	3.58\mbox{e-8}&	3.00\\
			\hline
		\end{tabular}
		\caption{$p=2$}
	\end{subtable}
	\begin{subtable}{\textwidth}
		\centering
		\begin{tabular}{||c||c|c||c|c||c|c||}
			\hline
			&\multicolumn{2}{c||}{CG}&\multicolumn{2}{c||}{VMS}&\multicolumn{2}{c||}{WENO}\\
			\hline
			$N_h$ & $\|u_h-u_{\text{exact}}\|_{L^1}$ & EOC& $\|u_h-u_{\text{exact}}\|_{L^1}$ & EOC & $\|u_h-u_{\text{exact}}\|_{L^1}$ & EOC\\
			\hline
			48  & 	3.46\mbox{e-4}&	  2.25&	2.53\mbox{e-4}&	 2.17&	2.29\mbox{e-4}&	1.91\\
			96  & 	1.42\mbox{e-5}&   4.60&	1.15\mbox{e-5}&	 4.46&	1.67\mbox{e-5}&	3.78\\
			192 & 	5.35\mbox{e-7}&   4.73&	4.92\mbox{e-7}&	 4.55&	1.00\mbox{e-6}&	4.06\\
			384 & 	2.77\mbox{e-8}&   4.27&	2.76\mbox{e-8}&	 4.16&	5.79\mbox{e-8}&	4.11\\
			768 &  	1.68\mbox{e-9}&   4.04&	1.67\mbox{e-9}&	 4.04&	3.49\mbox{e-9}&	4.05\\
			1536& 	1.05\mbox{e-10}&  4.01&	1.04\mbox{e-10}& 4.01&	2.15\mbox{e-10}&	4.02\\
			\hline
		\end{tabular}
		\caption{$p=3$}
	\end{subtable}
	\begin{subtable}{\textwidth}
		\centering
		\begin{tabular}{||c||c|c||c|c||c|c||}
			\hline
			&\multicolumn{2}{c||}{CG}&\multicolumn{2}{c||}{VMS}&\multicolumn{2}{c||}{WENO}\\
			\hline
			$N_h$ & $\|u_h-u_{\text{exact}}\|_{L^1}$ & EOC & $\|u_h-u_{\text{exact}}\|_{L^1}$ & EOC & $\|u_h-u_{\text{exact}}\|_{L^1}$ & EOC\\
			\hline
			64  &   5.07\mbox{e-5}&   3.98&   3.69\mbox{e-5}&   4.03&   1.15\mbox{e-4}&  3.78\\
			128 &   9.87\mbox{e-7}&   5.68&   4.83\mbox{e-7}&   6.26&   1.55\mbox{e-6}&  6.21\\
			256 &   9.11\mbox{e-8}&   3.44&   4.25\mbox{e-8}&   3.51&   6.19\mbox{e-8}&  4.65\\
			512 &   5.47\mbox{e-9}&   4.06&   1.59\mbox{e-9}&   4.74&	1.27\mbox{e-9}&  5.60\\
			1024&   3.32\mbox{e-10}&  4.04&   5.89\mbox{e-11}&  4.76&	3.70\mbox{e-11}&  5.10\\
			\hline
		\end{tabular}
		\caption{$p=4$}
	\end{subtable}
	\caption{1D Burgers equation, grid convergence history for finite elements of degree $p\in\{1,2,3,4\}$.}
	\label{tab:convbur}
\end{table}

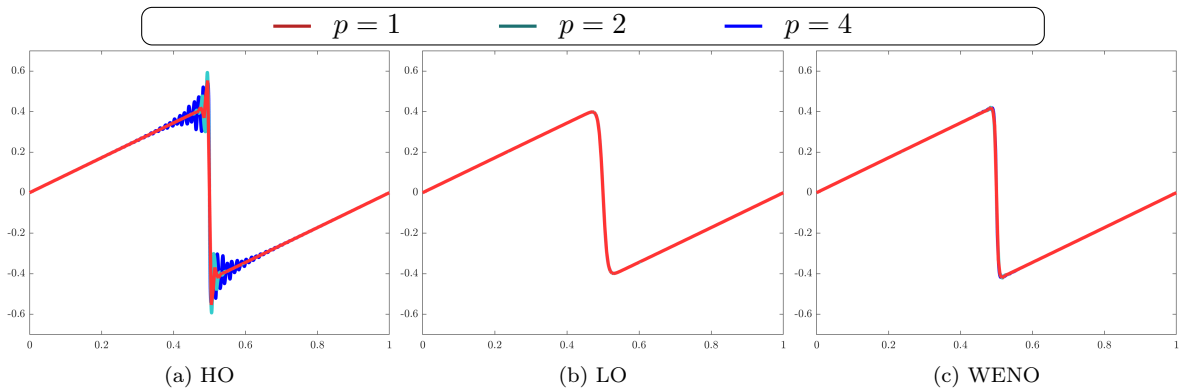
\begin{figure}[!htb]
	\centering
	\begin{subfigure}[b]{.9\linewidth}
		\centering
		\begin{tikzpicture}
		\draw[rounded corners] (0, 0) rectangle (12, 0.5) node[pos=.5]{};
		\draw[very thick, color={rgb:red,1;green,0.2;blue,0.2}] (1.75,0.25)--(2.25,0.25);
		\node at (3,0.25) (a) {$p=1$};
		\node at (6,0.25) (a) {$p=2$};
		\node at (9,0.25) (a) {$p=4$};
		\draw[very thick,color={rgb:red,0.2;green,0.8;blue,0.8}] (4.75,0.25)--(5.25,0.25);
		\draw[very thick,color={rgb:red,0;green,0;blue,1}] (7.75,0.25)--(8.25,0.25);
		\end{tikzpicture}
	\end{subfigure}
	\begin{subfigure}[b]{.3\linewidth}
		\includegraphics[width=\linewidth]{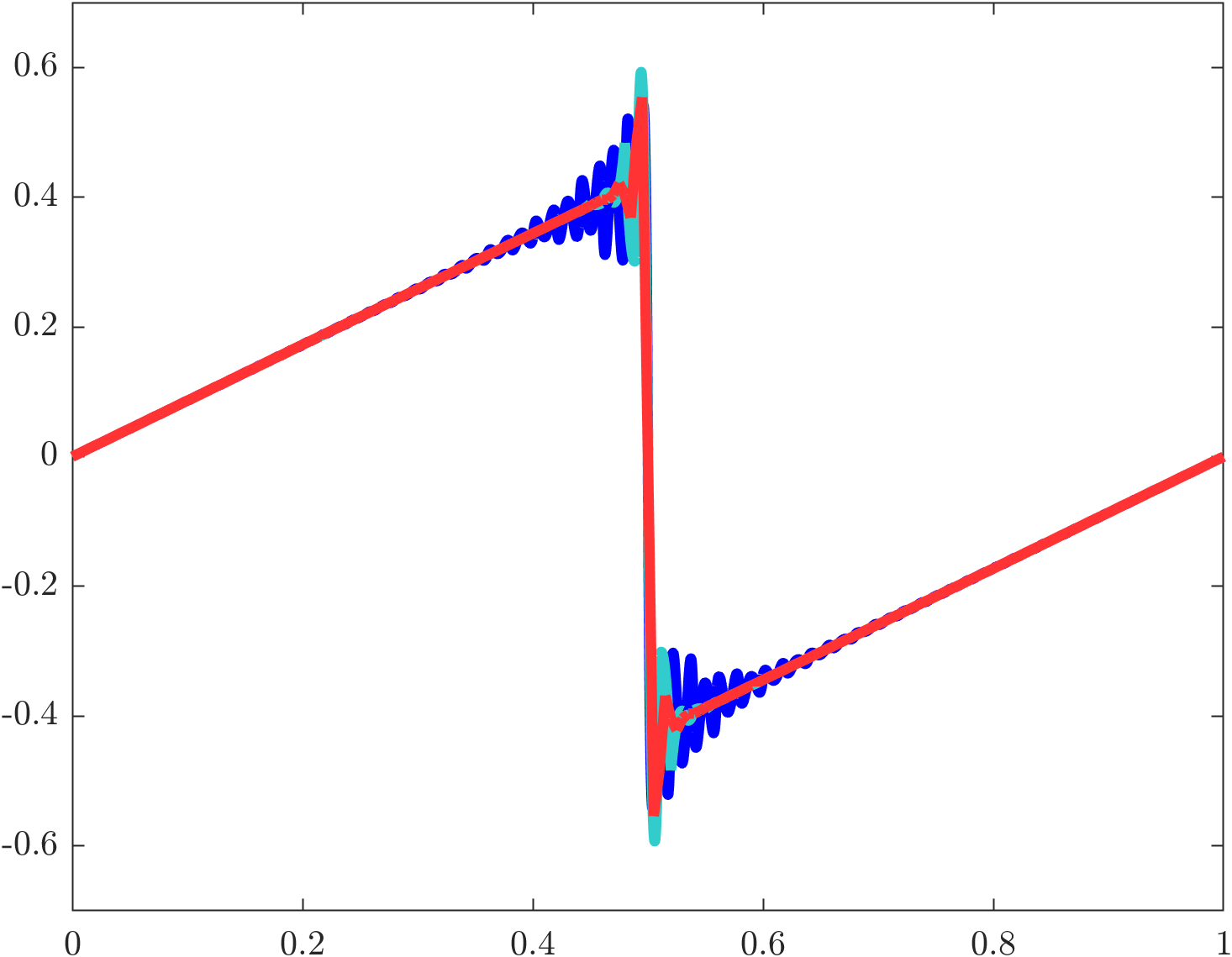}
		\caption{HO}
		\label{fig:burho}
	\end{subfigure}
	\begin{subfigure}[b]{.3\linewidth}
		\includegraphics[width=\linewidth]{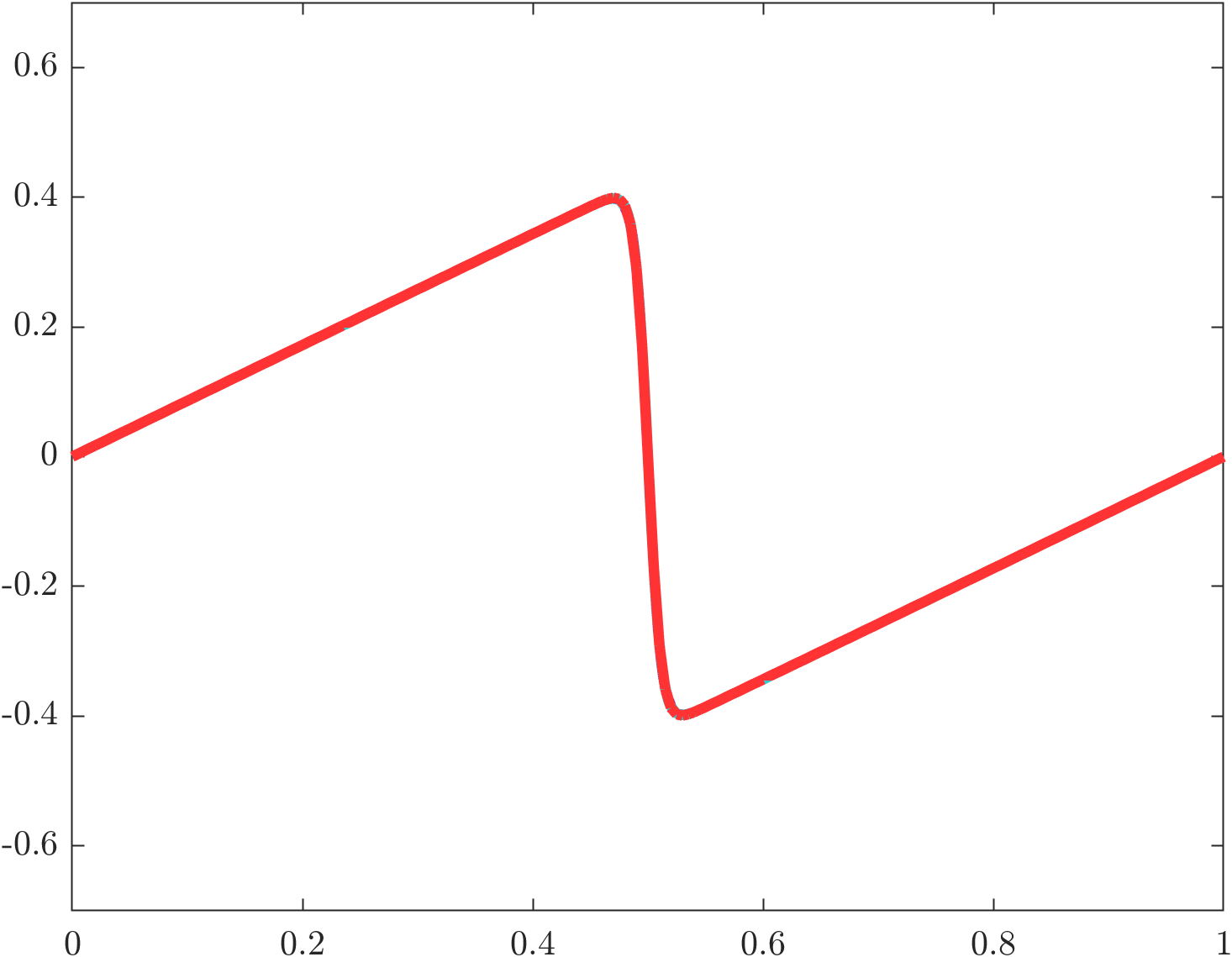}
		\caption{LO}
		\label{fig:burlo}
	\end{subfigure}
	\begin{subfigure}[b]{.3\linewidth}
		\includegraphics[width=\linewidth]{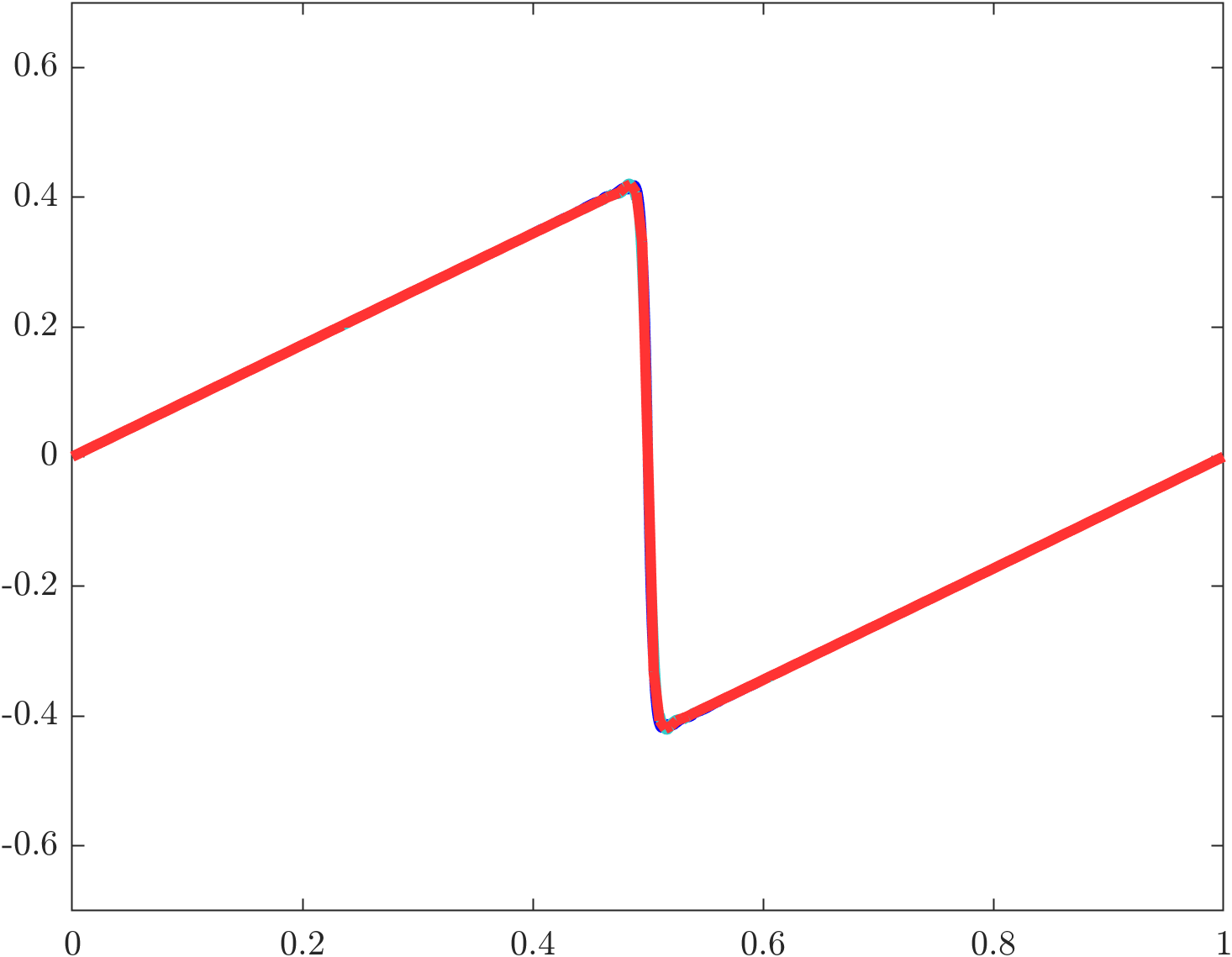}
		\caption{WENO}
		\label{fig:burweno}
	\end{subfigure}
	\caption{1D Burgers equation, numerical solutions at $t=1.0$ obtained using $N_h=200$ and $p=\{1,2,4\}$.}
	\label{fig:burgers}
\end{figure}

Extending the final time to $t=1.0>t_c$, we study the shock-capturing properties of the WENO scheme and its linear ingredients. Simulations are run using $p\in \{1,2,4\}$ on meshes corresponding to $N_h=200$ DoFs. The HO results shown in Fig. \ref{fig:burho} indicate that high-order dissipation is not suited for shock capturing. The LO solutions presented in Fig. \ref{fig:burlo} are nonoscillatory and not as diffusive as in the case of linear advection, because shock fronts are self-steepening. Even higher resolution can be achieved using the WENO scheme, which produces the result shown in Fig. \ref{fig:burweno}. Remarkably, the LO and WENO approximations remain virtually unchanged as we vary $p$ while keeping $N_h$ fixed. The differences between the oscillatory HO results for different polynomial degrees are clearly visible.

\subsection{Two-dimensional solid body rotation}

LeVeque's \cite{leveque1996} solid body rotation problem is a popular stability test for discretizations of
\begin{align}
\frac{\partial u}{\partial t}+\nabla\cdot(\mathbf{v}u) = 0 \quad \mbox{in}\ \Omega=(0,1)^2.
\label{num:sbreq}
\end{align}
The velocity field ${\bf v}(x,y)=2\pi(0.5-y,x-0.5)$ is divergence-free. The initial
condition is given by
\begin{align}
u_0(x,y) = \begin{cases}
u_0^{\text{hump}}(x,y) & \mbox{if} \ \sqrt{(x-0.25)^2+(y-0.5)^2}\le0.15, \\
u_0^{\text{cone}}(x,y) & \mbox{if} \ \sqrt{(x-0.5)^2+(y-0.25)^2}\le0.15, \\
1 & \mbox{if} \; \begin{cases}
\sqrt{(x-0.5)^2+(y-0.75)^2}\le0.15, \\
|x-0.5|\ge0.025,y\ge0.85,
\end{cases}\\
0 & \mbox{otherwise},
\end{cases}
\label{num:sbrinit}
\end{align}
where
\begin{align*}
u_0^{\text{hump}}(x,y) &= \frac{1}{4}+\frac{1}{4}\cos\bigg(\frac{\pi\sqrt{(x-0.25)^2+(y-0.5)^2}}{0.15}\bigg), \\
u_0^{\text{cone}}(x,y) &= 1-\frac{\sqrt{(x-0.5)^2+(y-0.25)^2}}{0.15}.
\end{align*}
The initial configuration rotates around the center
$(0.5,0.5)$ of $\Omega$.
After each complete revolution (i.e., for $t\in\mathbb{N}$), the exact solution
of the linear advection equation \eqref{num:bureq}
coincides with $u_0$.

\begin{figure}[!htb]
	\centering
	\begin{subfigure}[b]{.45\linewidth}
		\includegraphics[width=\linewidth]{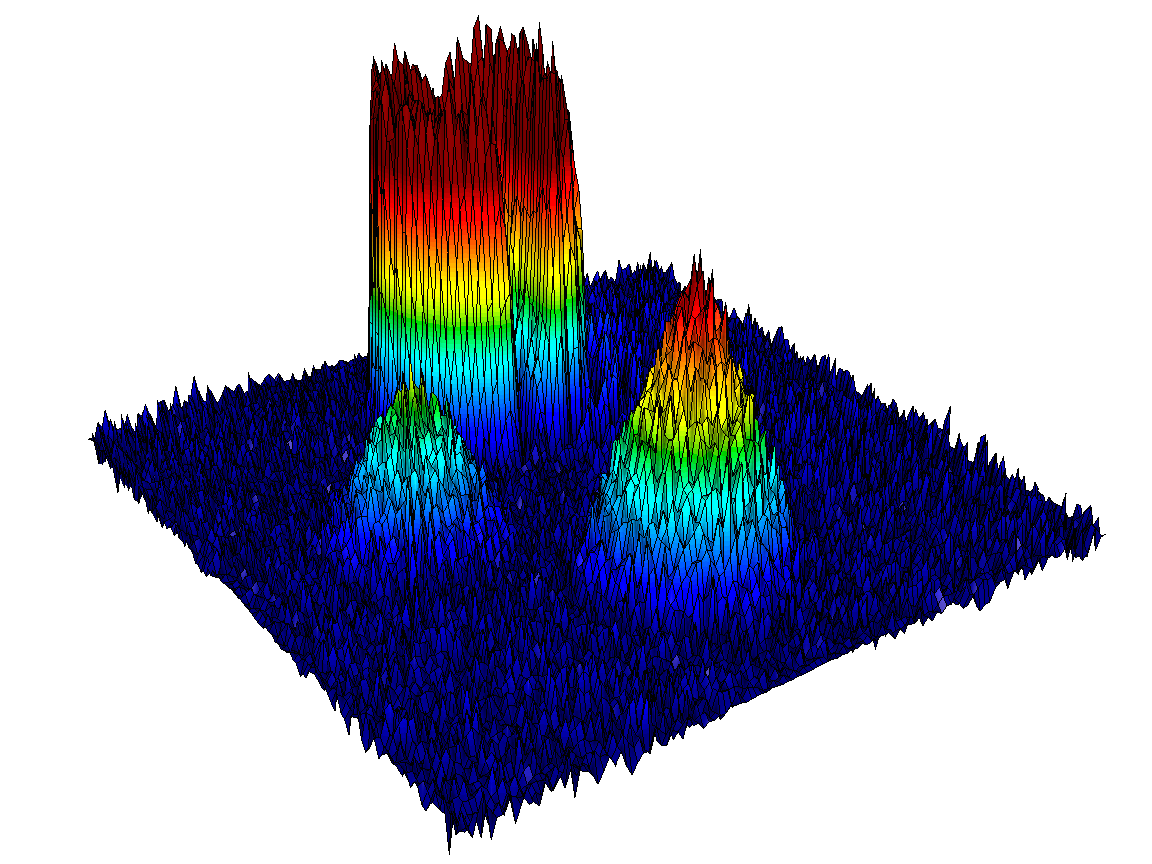}
		\caption{CG, $E_1$=3.70e-02, $u_h\in [-0.399,1.367]$}
	\end{subfigure}
	\begin{subfigure}[b]{.45\linewidth}
		\includegraphics[width=\linewidth]{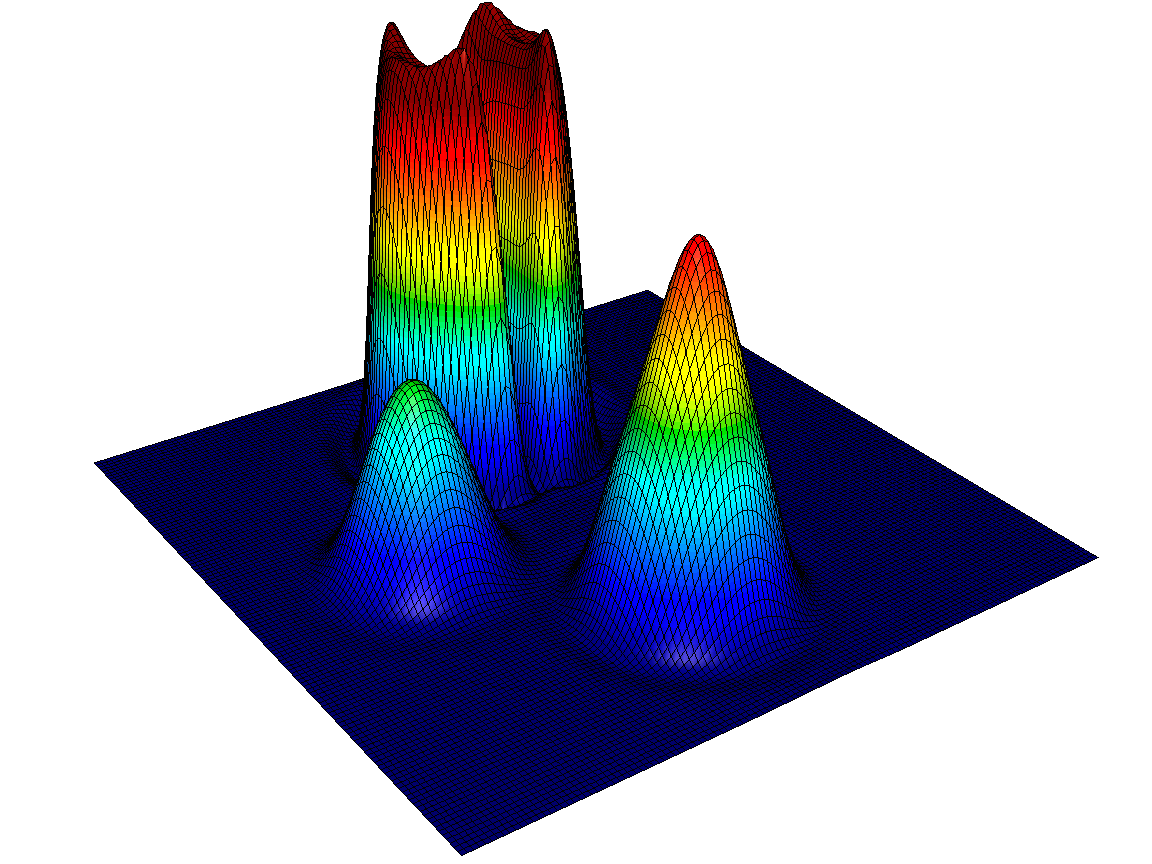}
		\caption{HO, $E_1$=1.47e-02, $u_h \in [-0.055,1.122]$}
	\end{subfigure}
	\begin{subfigure}[b]{.45\linewidth}
		\includegraphics[width=\linewidth]{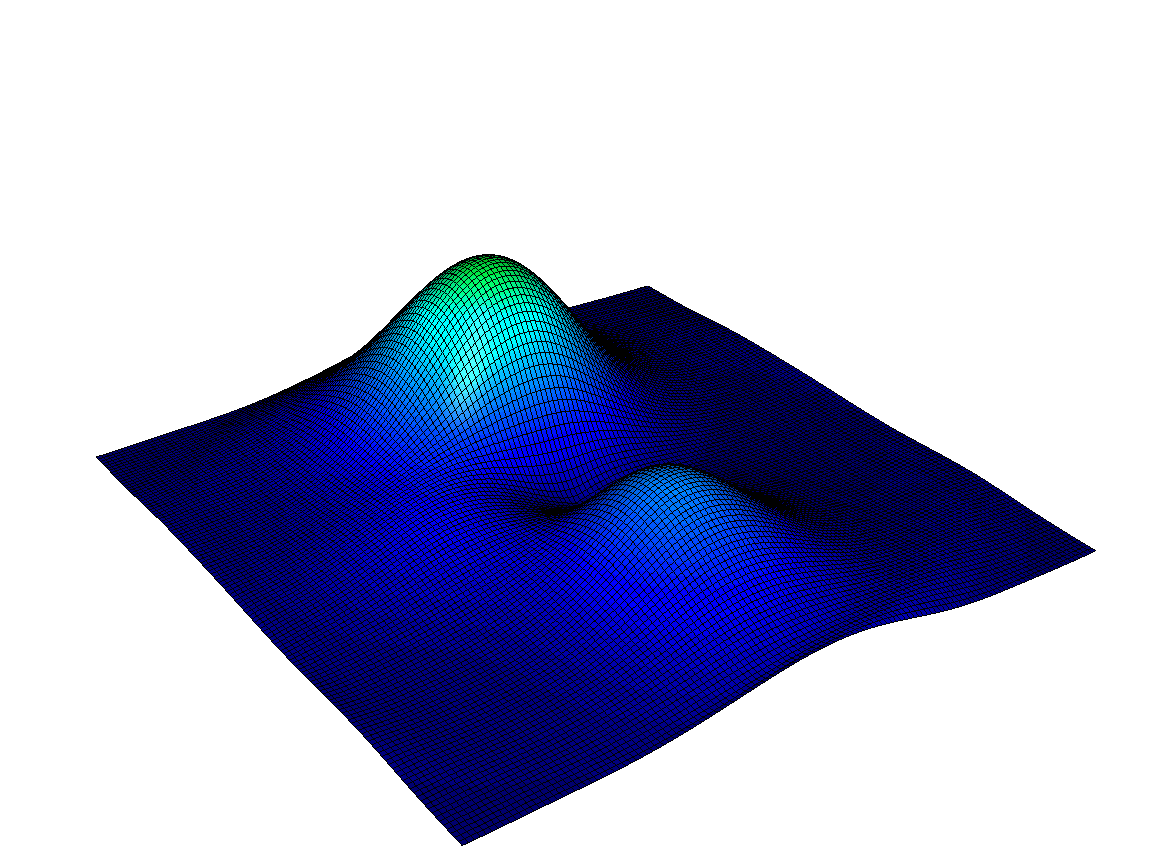}
		\caption{LO, $E_1$=1.10e-01, $u_h \in [0.003,0.478]$}
	\end{subfigure}
	\begin{subfigure}[b]{.45\linewidth}
		\includegraphics[width=\linewidth]{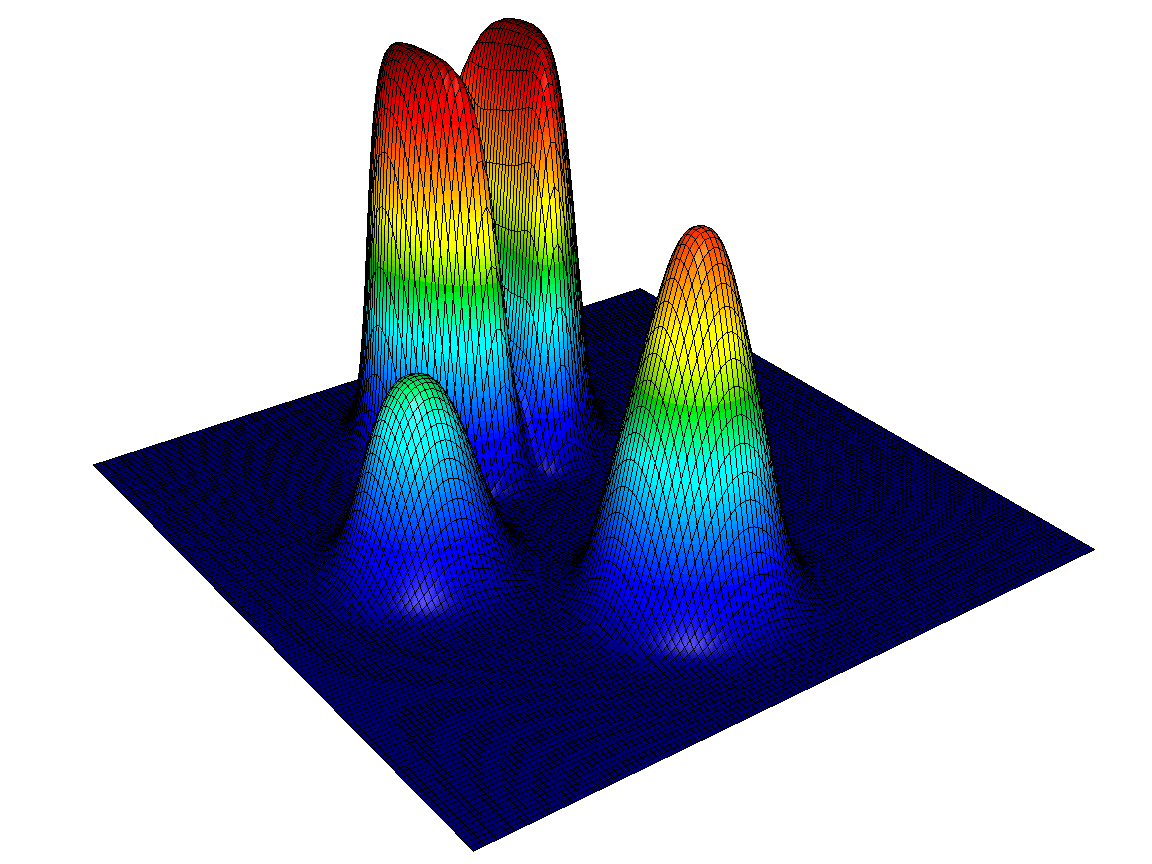}
		\caption{WENO, $E_1$=2.68e-02, $u_h \in [-0.001,0.993]$}
	\end{subfigure}
	\caption{Solid body rotation, numerical solutions at $t=1$ obtained using $N_h=129^2$ and $p=1$.}
	\label{fig:sbrp1}
\end{figure}

\begin{figure}[!htb]
	\centering
	\begin{subfigure}[b]{.45\linewidth}
		\includegraphics[width=\linewidth]{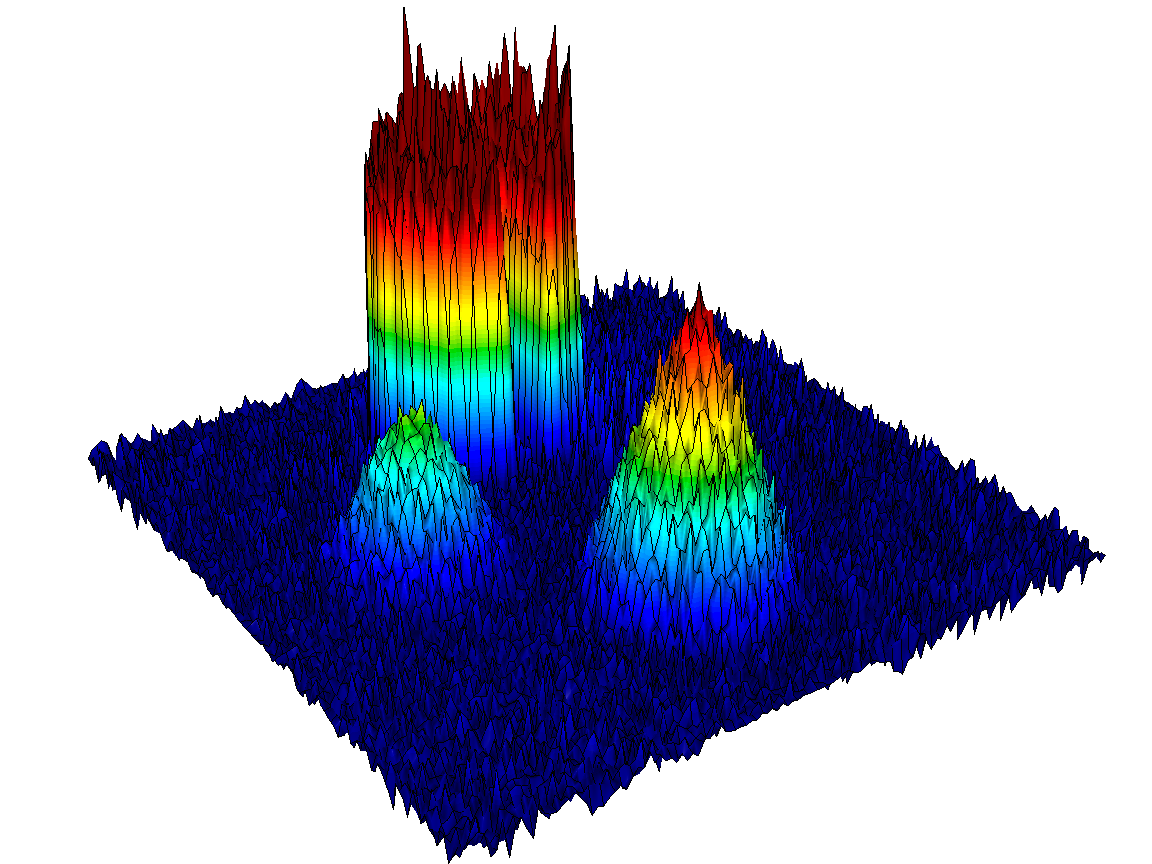}
		\caption{CG, $E_1$=3.17e-02, $u_h\in [-0.377,1.476]$}
	\end{subfigure}
	\begin{subfigure}[b]{.45\linewidth}
		\includegraphics[width=\linewidth]{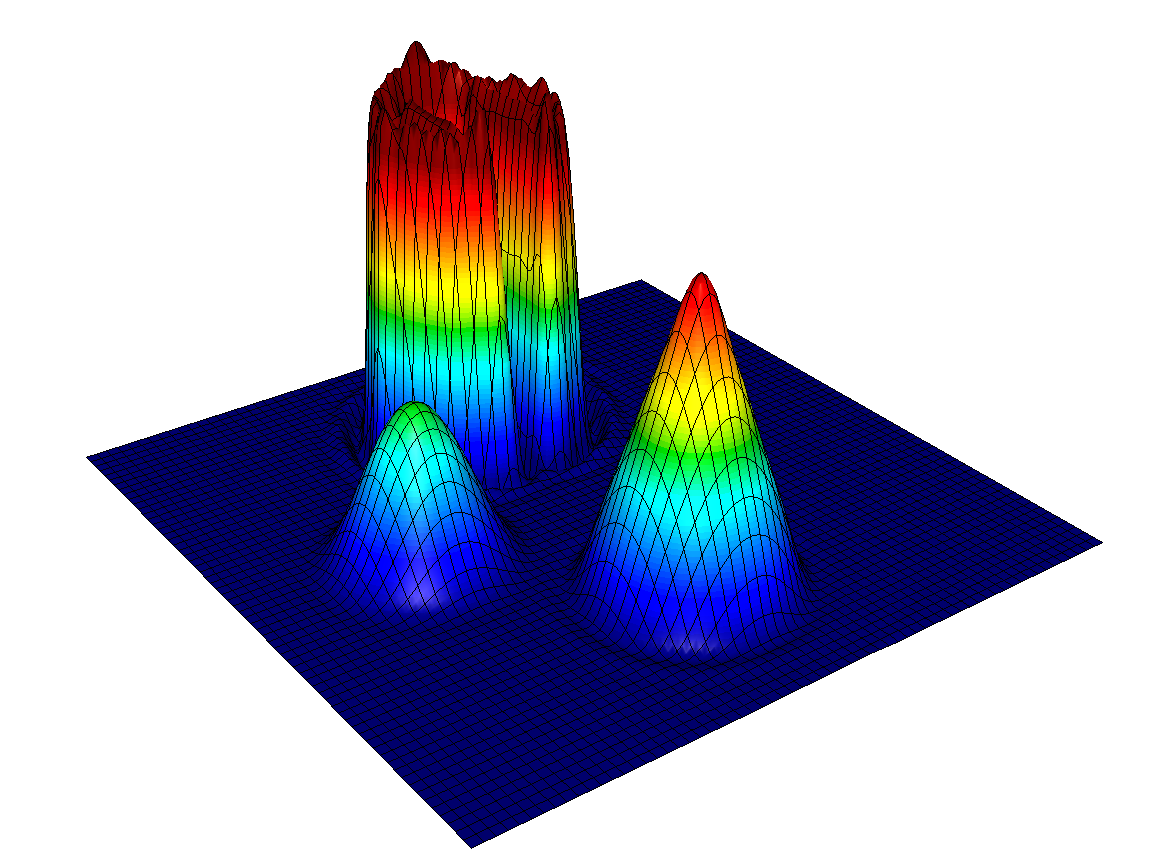}
		\caption{HO, $E_1$=1.23e-02, $u_h\in [-0.213,1.183]$}
	\end{subfigure}
	\begin{subfigure}[b]{.45\linewidth}
		\includegraphics[width=\linewidth]{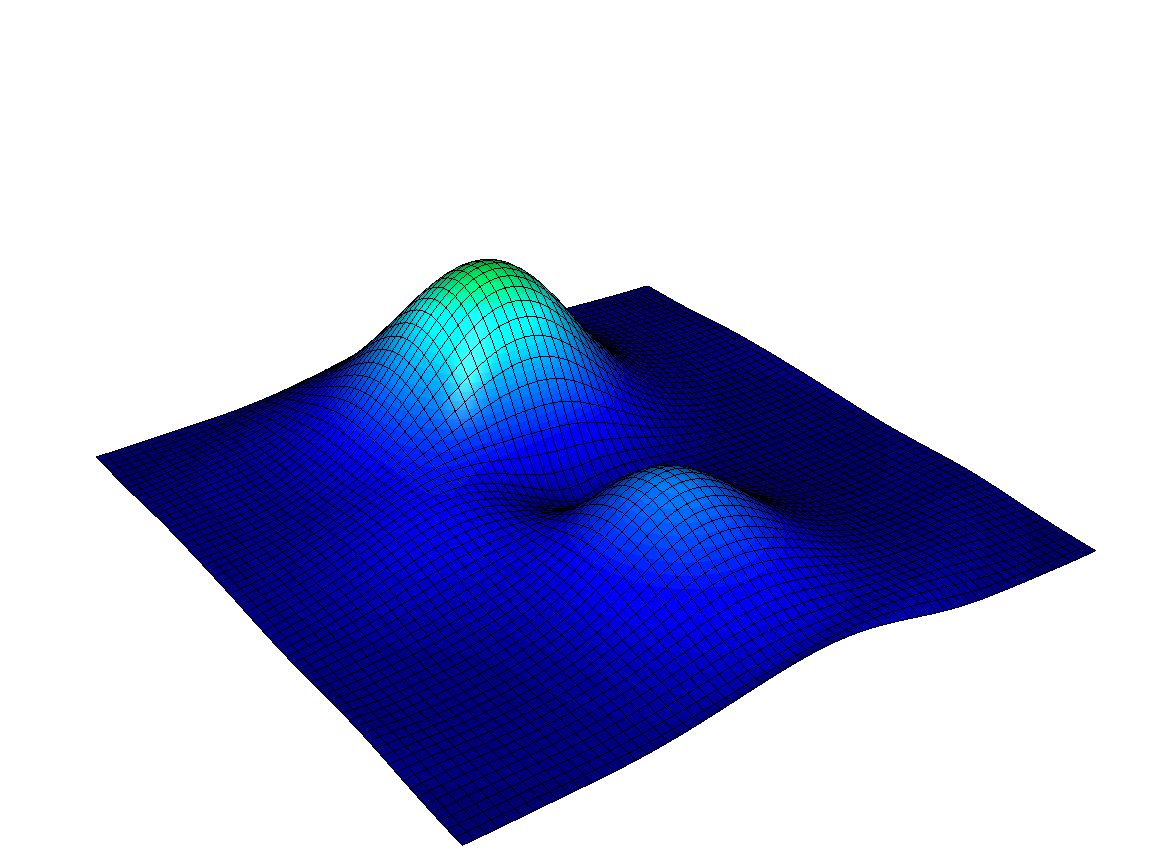}
		\caption{LO, $E_1$=1.10e-01, $u_h \in [0.004,0.465]$}
	\end{subfigure}
	\begin{subfigure}[b]{.45\linewidth}
		\includegraphics[width=\linewidth]{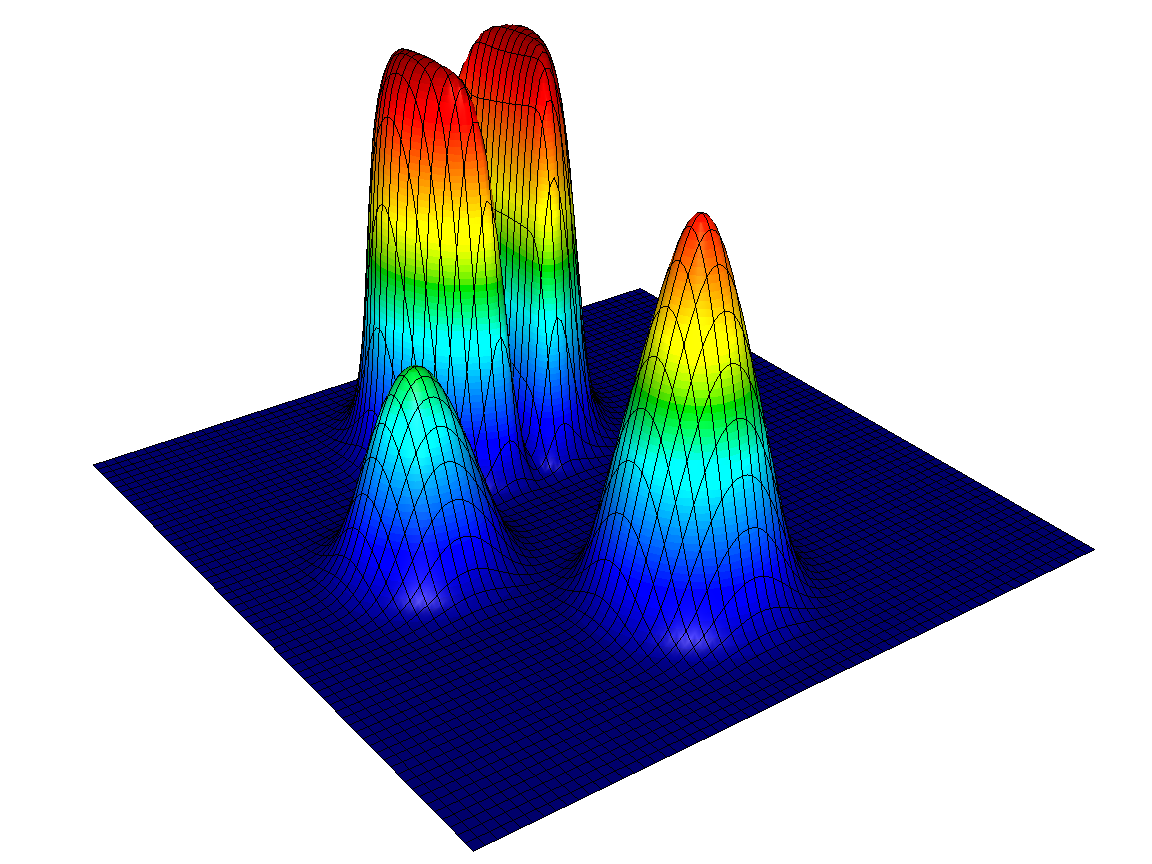}
		\caption{WENO, $E_1$=2.67e-01, $u_h \in [0.000,0.981]$}
	\end{subfigure}
	\caption{Solid body rotation, numerical solutions at $t=1$ obtained using $N_h=129^2$ and $p=2$.}
	\label{fig:sbrp2}
\end{figure}

\begin{figure}[!htb]
	\centering
	\begin{subfigure}[b]{.45\linewidth}
		\includegraphics[width=\linewidth]{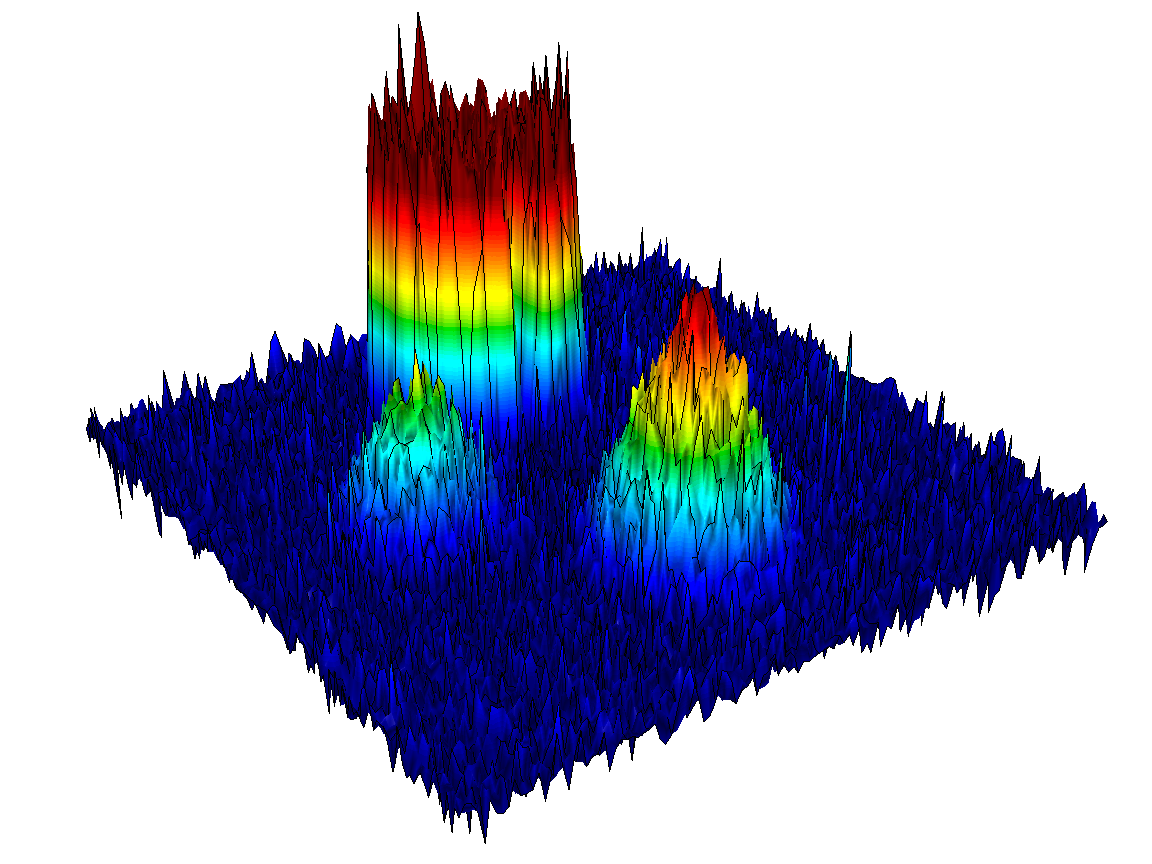}
		\caption{CG, $E_1$=4.19e-02, $u_h \in [-0.480,1.429]$}
	\end{subfigure}
	\begin{subfigure}[b]{.45\linewidth}
		\includegraphics[width=\linewidth]{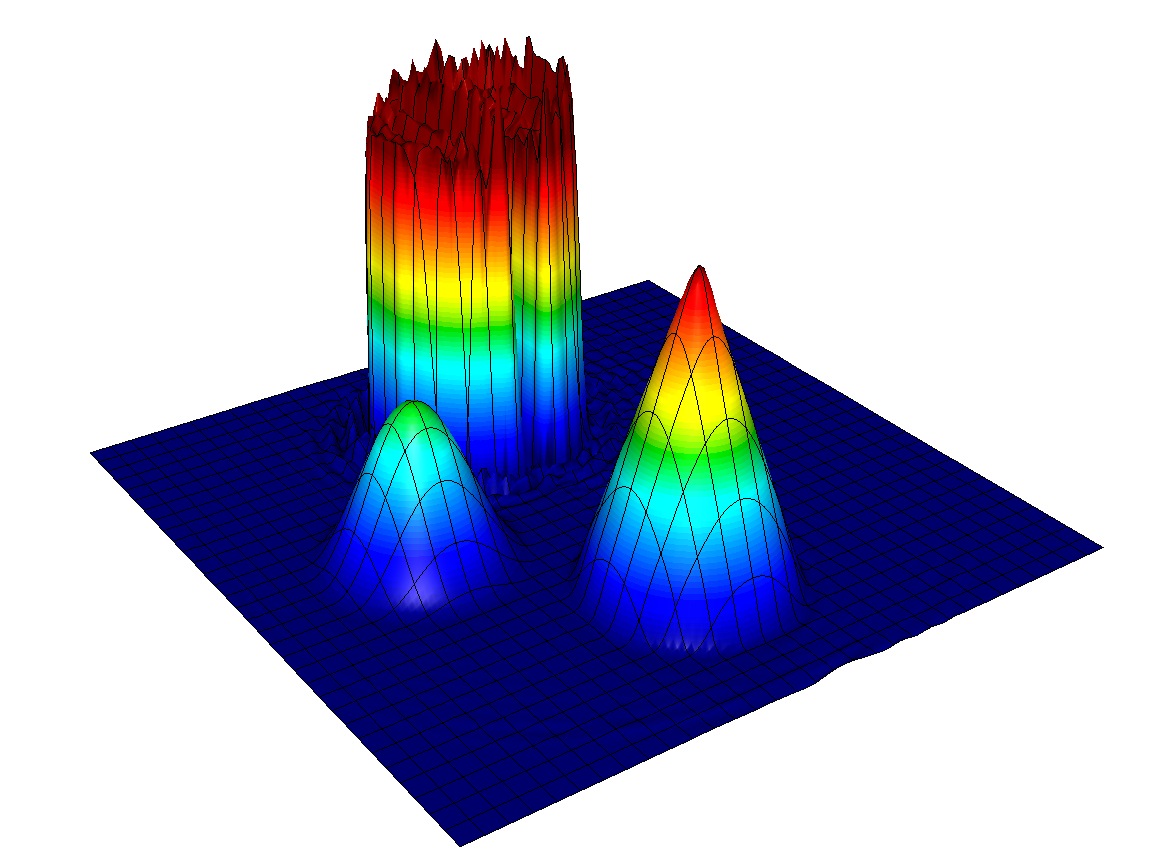}
		\caption{HO, $E_1$=1.25e-02, $u_h \in [-0.265,1.266]$}
	\end{subfigure}
	\begin{subfigure}[b]{.45\linewidth}
		\includegraphics[width=\linewidth]{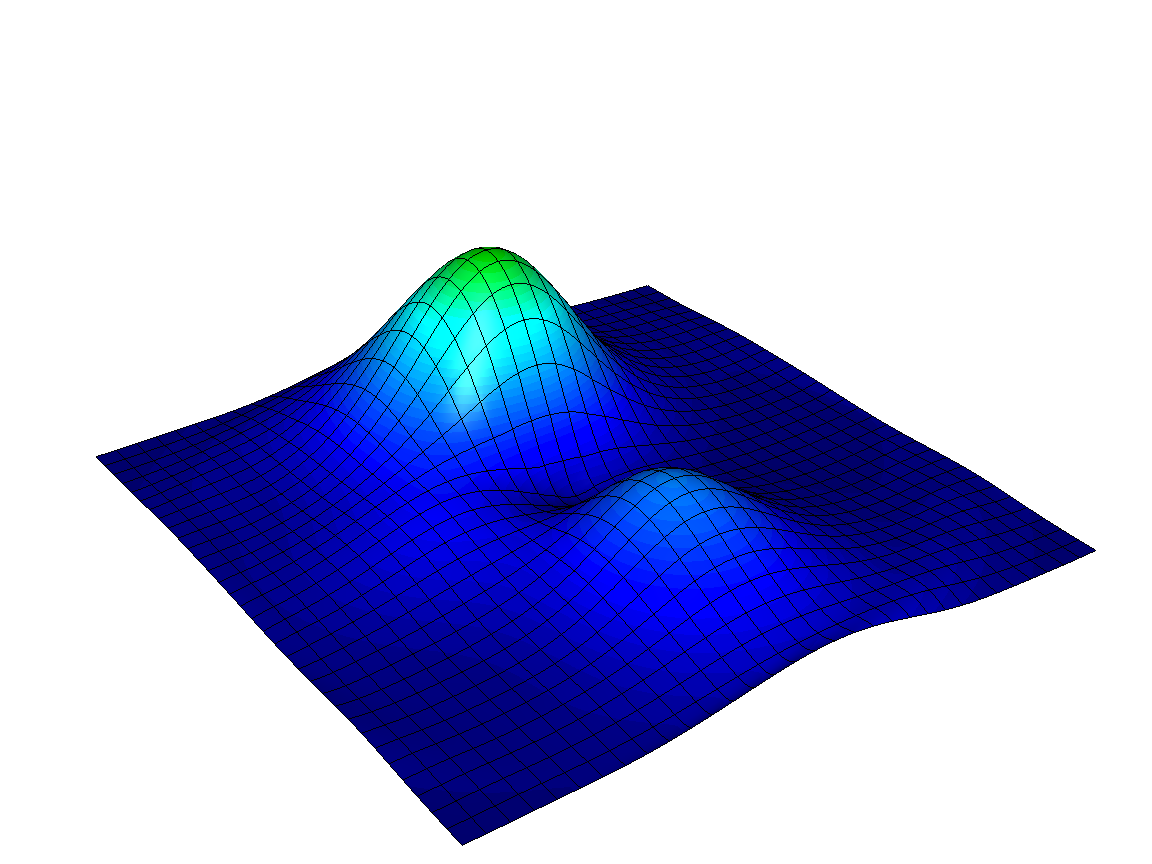}
		\caption{LO, $E_1$=1.11e-01, $u_h \in [0.004,0.501]$}
	\end{subfigure}
	\begin{subfigure}[b]{.45\linewidth}
		\includegraphics[width=\linewidth]{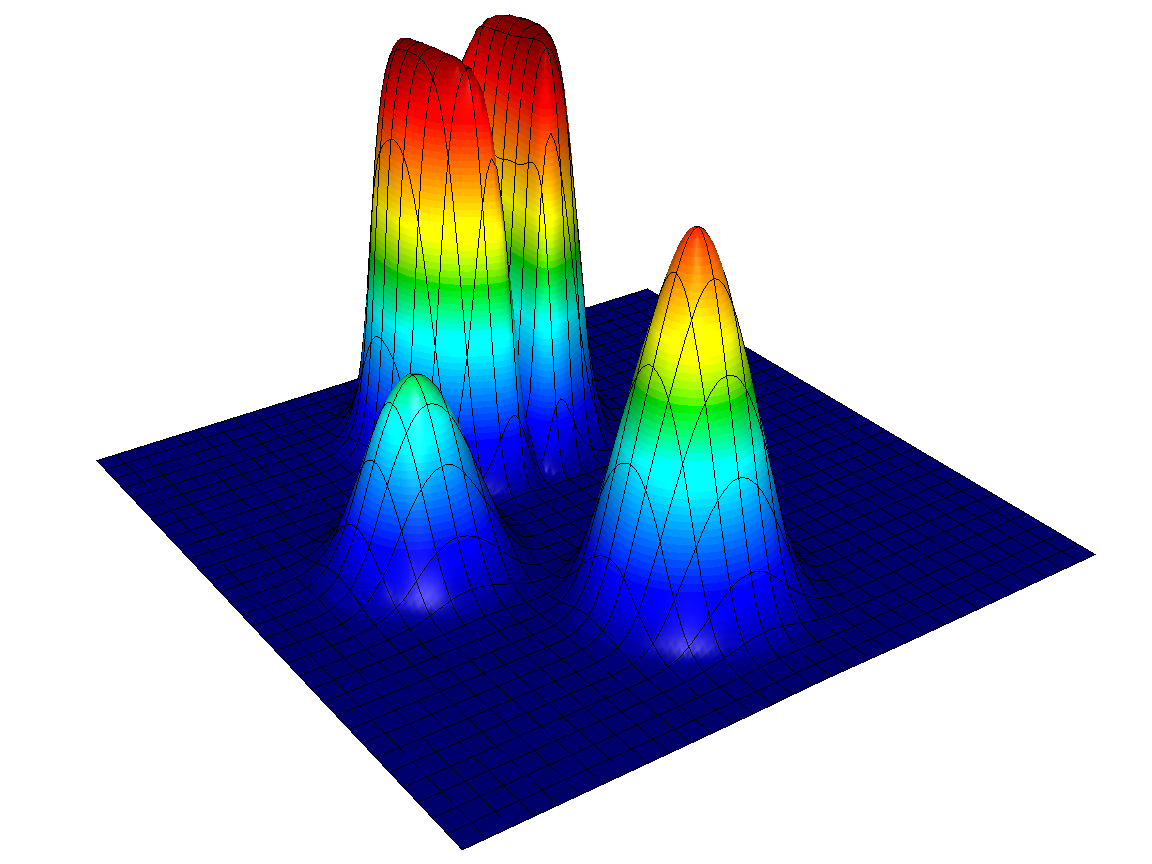}
		\caption{WENO, $E_1$=2.76e-02, $u_h \in [0.000,0.998]$}
	\end{subfigure}
	\caption{Solid body rotation, numerical solutions at $t=1$ obtained using $N_h=129^2$ and $p=4$.}
	\label{fig:sbrp4}
\end{figure}

We evolve numerical solutions up to the finite time $t=1.0$ on uniform Cartesian meshes
corresponding to $N_h=129^2$ DoFs. The results for Lagrange finite elements of degree $p=\{1,2,4\}$
are presented in  \mbox{Figs~\ref{fig:sbrp1}--\ref{fig:sbrp4}}. All
CG approximations are corrupted by global spurious oscillations. The HO stabilization is
sufficient in smooth regions but large undershoots and overshoots are observed in the vicinity of the
slotted cylinder. The LO solution is nonoscillatory but very diffusive. Using the
smoothness indicator $\gamma_e$ to adaptively blend HO and LO dissipation, the nonlinear
WENO scheme exploits the complementary advantages of the two linear methods and
yields superior results for all combinations of the mesh size $h$ and polynomial
degree $p$. To show that $h$ is coarsened as $p$ is refined to keep $N_h$ fixed,
we plot the edges of the meshes on which the computations were performed.

\subsection{KPP problem}
Finally, we consider the KPP problem \cite{kurganov2007,kuzmin2020g}, a challenging nonlinear test for assessment of entropy stability properties. The conservation law \eqref{ibvp-pde} with 
the nonconvex flux function
\beq\label{num:kppeq}
\mathbf{f}(u)=(\sin(u),\cos(u))
\eeq
is solved in the computational domain
$\Omega_h=(-2,2)\times(-2.5,1.5)$ using the initial
condition
\beq \label{num:kppinit}
u_0(x,y)=\begin{cases}
\frac{7\pi}{2} & \mbox{if}\quad \sqrt{x^2+y^2}\le 1,\\
\frac{\pi}{4} & \mbox{otherwise}.
\end{cases}
\eeq
A global upper bound for the maximum speed that we need to calculate the viscosity parameter $\nu_e$ is given by $\lambda_e=1$. To make the LO component of the WENO scheme as dissipative as necessary to safely suppress undershoots/overshoots even on coarsest meshes, we define it using $\lambda_e=2$. Deviating from the default settings, we choose the linear weights
$\tilde \omega_l^e = 0.2,\ l=0,\ldots,4$. As mentioned above, smaller values
of $\tilde\omega_0^e$ make the smoothness indicator $\gamma_e$ more sensitive to spurious oscillations. Finally, we lump the mass matrix in the CG version without any stabilization ($\nu_e=0$) because the consistent-mass CG approximation was found to produce extremely oscillatory results in this test.

\begin{figure}[!htb]
	\centering
	\begin{subfigure}[b]{.45\linewidth}
		\includegraphics[width=\linewidth,trim=0 50 0 0,clip]{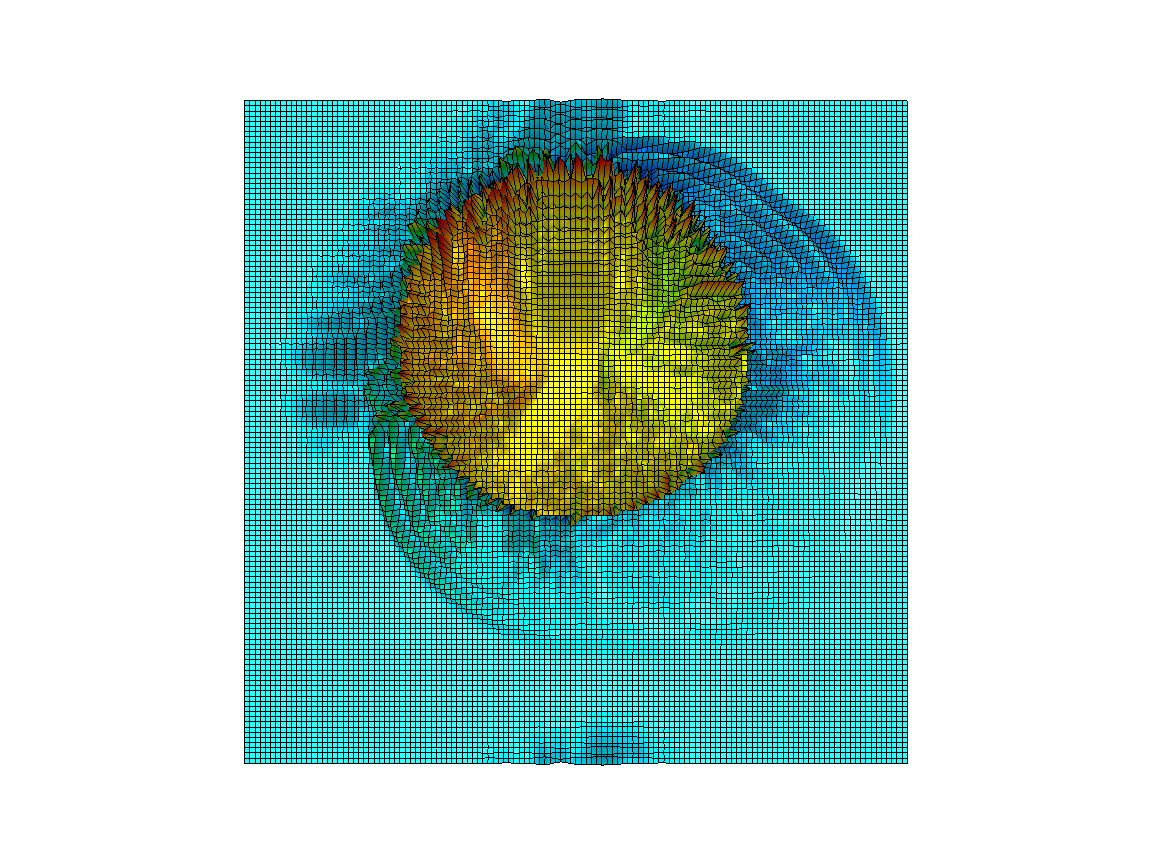}
		\caption{CG, $u_h \in [-11.752, 24.446]$}
	\end{subfigure}
	\begin{subfigure}[b]{.45\linewidth}
		\includegraphics[width=\linewidth,trim=0 50 0 0,clip]{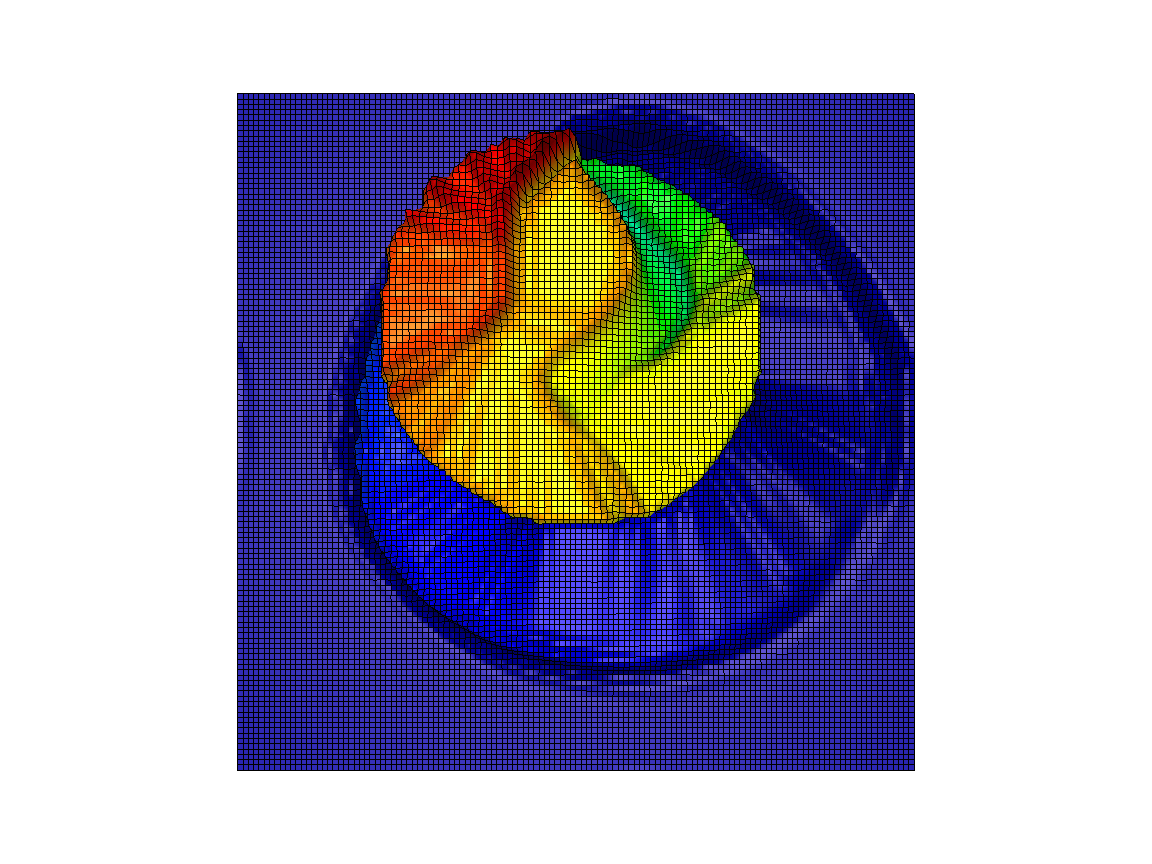}
		\caption{HO, $u_h \in [-2.986, 15.271]$}
	\end{subfigure}
	\begin{subfigure}[b]{.45\linewidth}
		\includegraphics[width=\linewidth,trim=0 50 0 0,clip]{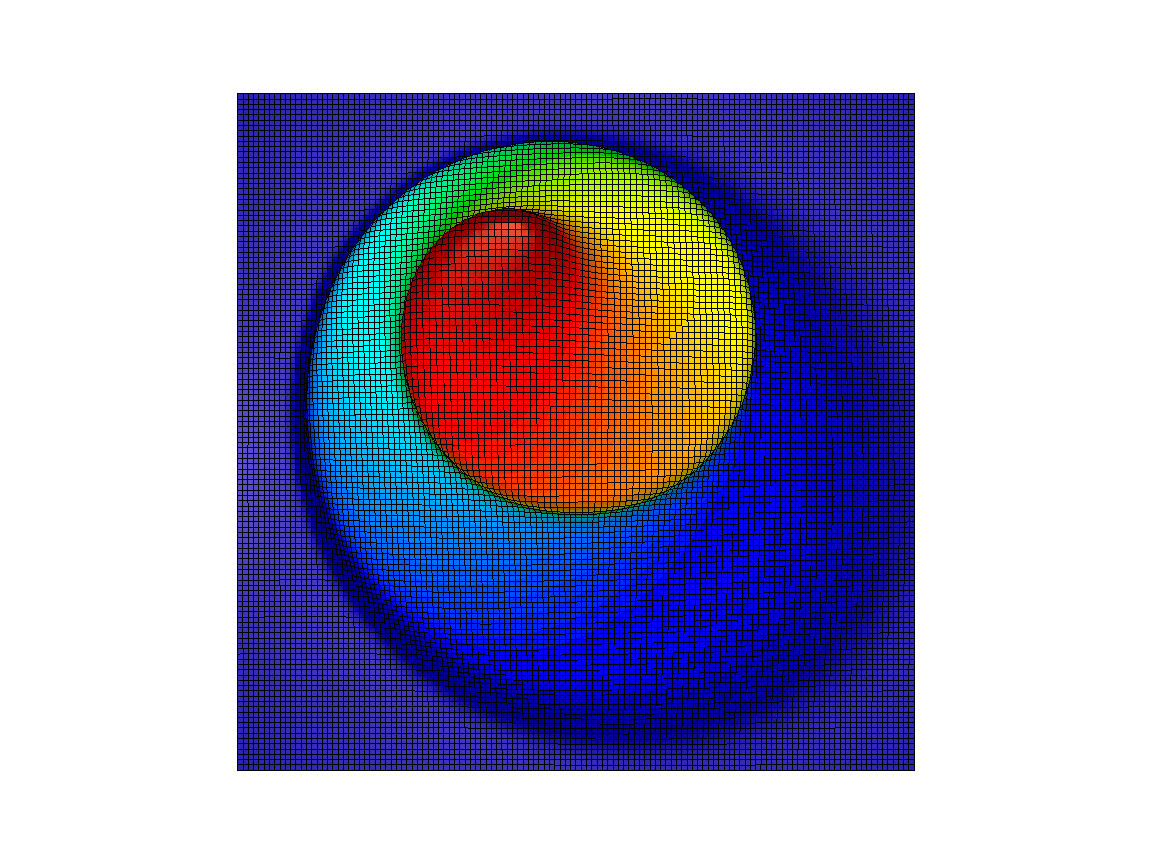}
		\caption{LO, $u_h \in [0.785,10.893]$}
	\end{subfigure}
	\begin{subfigure}[b]{.45\linewidth}
		\includegraphics[width=\linewidth,trim=0 50 0 0,clip]{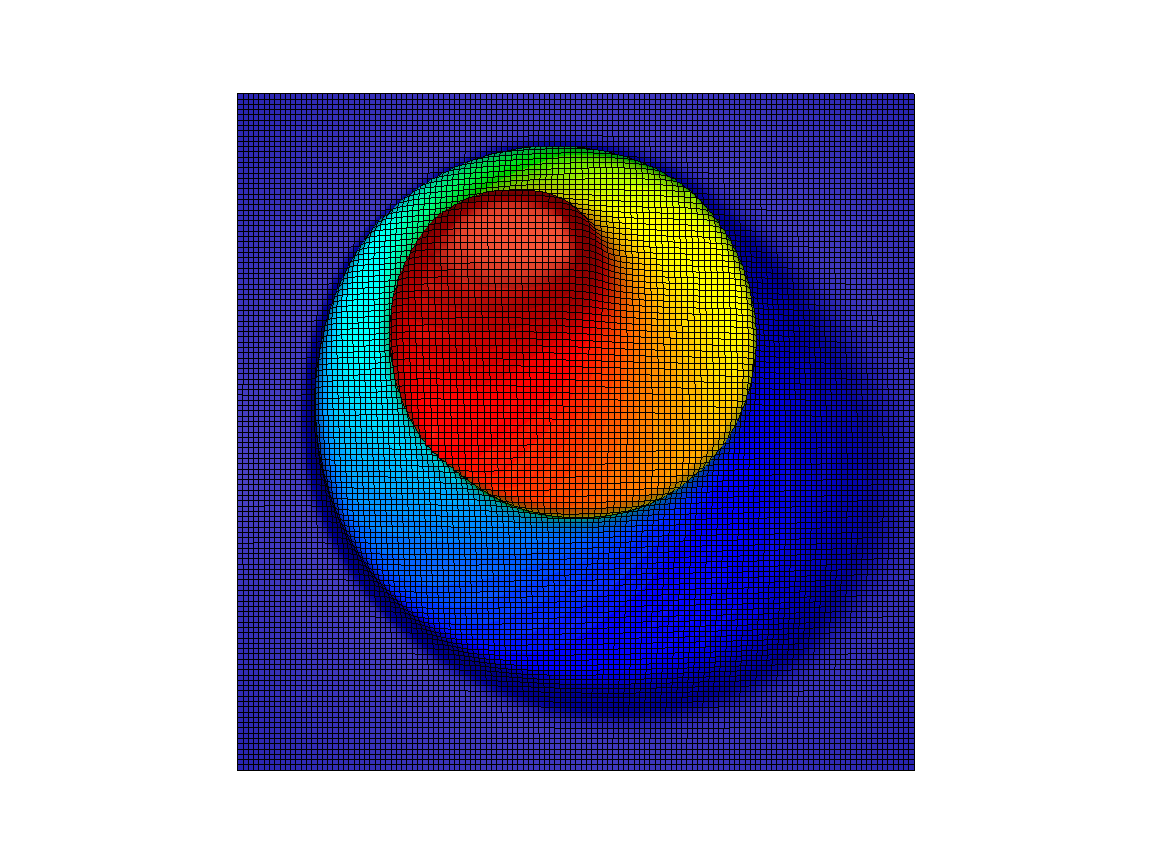}
		\caption{WENO, $u_h \in [0.785,10.990]$}
	\end{subfigure}
	\caption{KPP problem, numerical solutions at $t=1$ obtained using $N_h=129^2$ and $p=1$.}
	\label{fig:kppp1}
\end{figure}

\begin{figure}[!htb]
	\centering
	\begin{subfigure}[b]{.45\linewidth}
		\includegraphics[width=\linewidth,trim=0 50 0 0,clip]{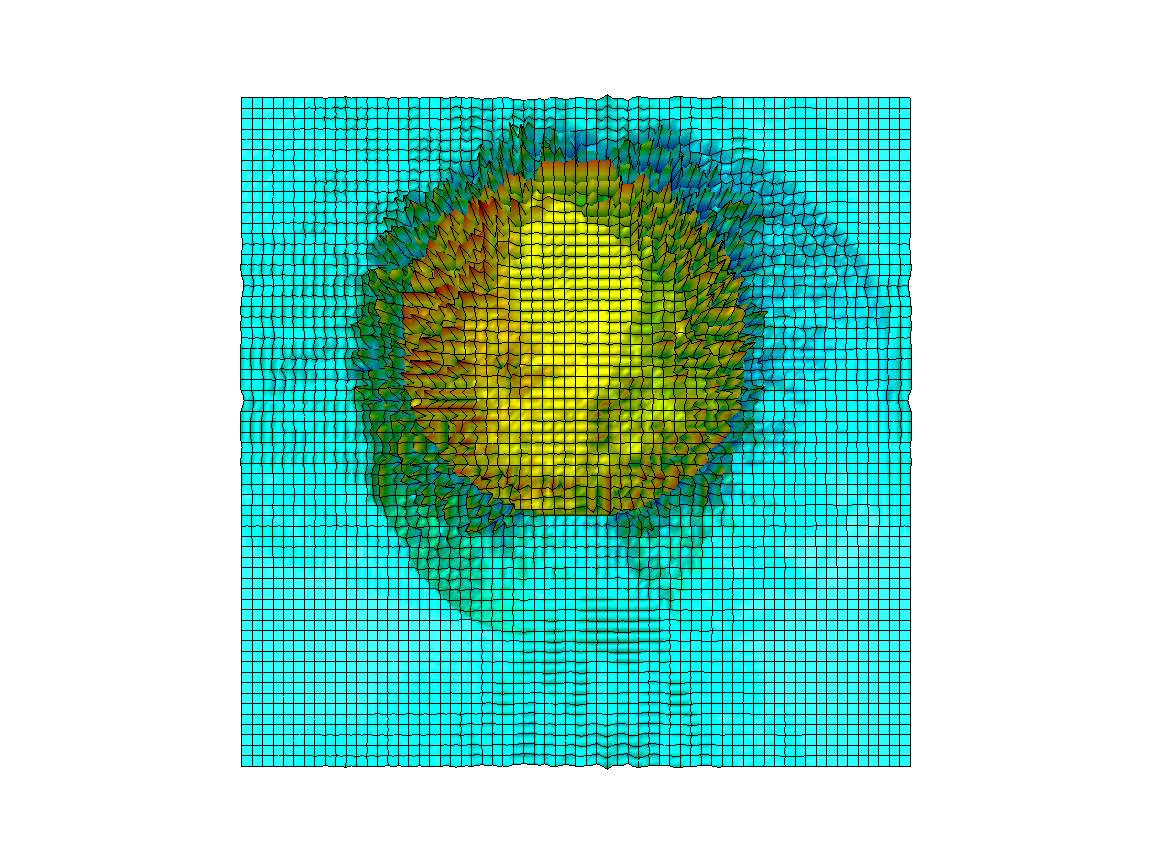}
		\caption{CG, $u_h \in [-16.934,29.237]$}
	\end{subfigure}
	\begin{subfigure}[b]{.45\linewidth}
		\includegraphics[width=\linewidth,trim=0 50 0 0,clip]{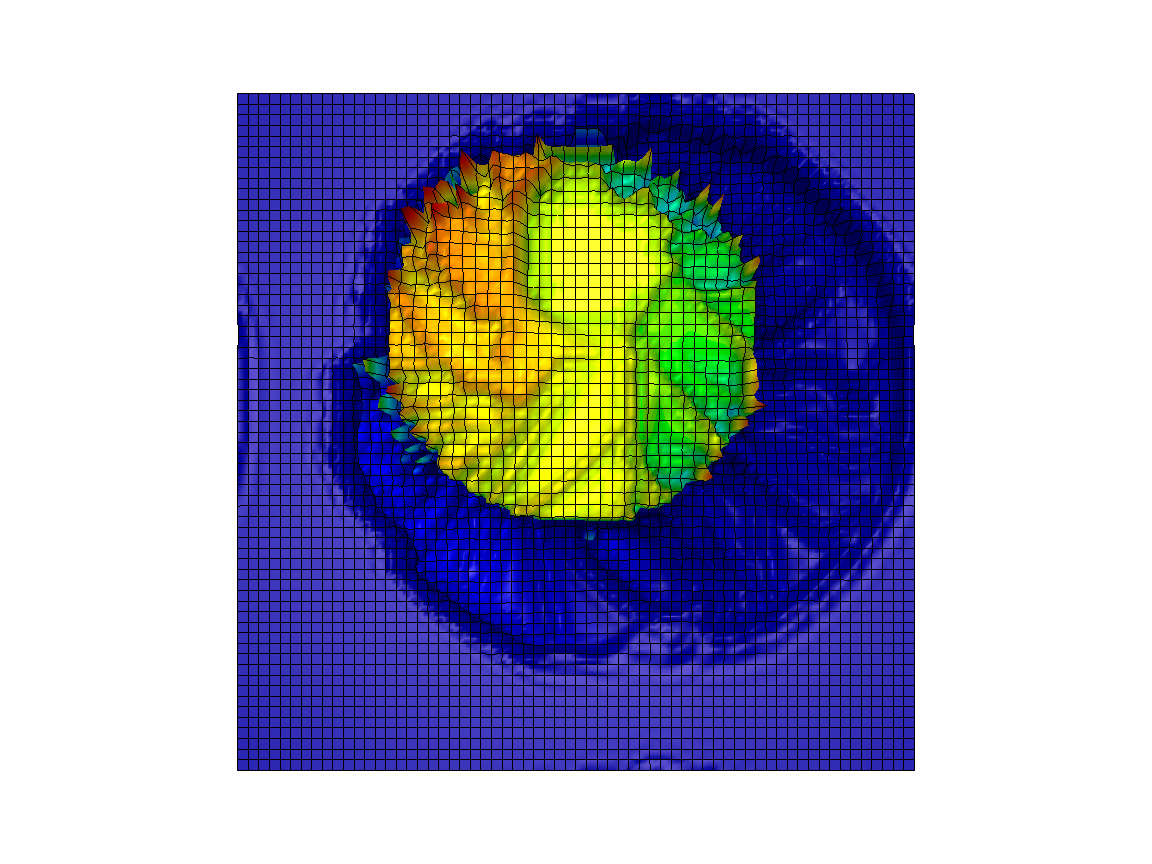}
		\caption{HO, $u_h \in [-8.503, 19.027]$}
	\end{subfigure}
	\begin{subfigure}[b]{.45\linewidth}
		\includegraphics[width=\linewidth,trim=0 50 0 0,clip]{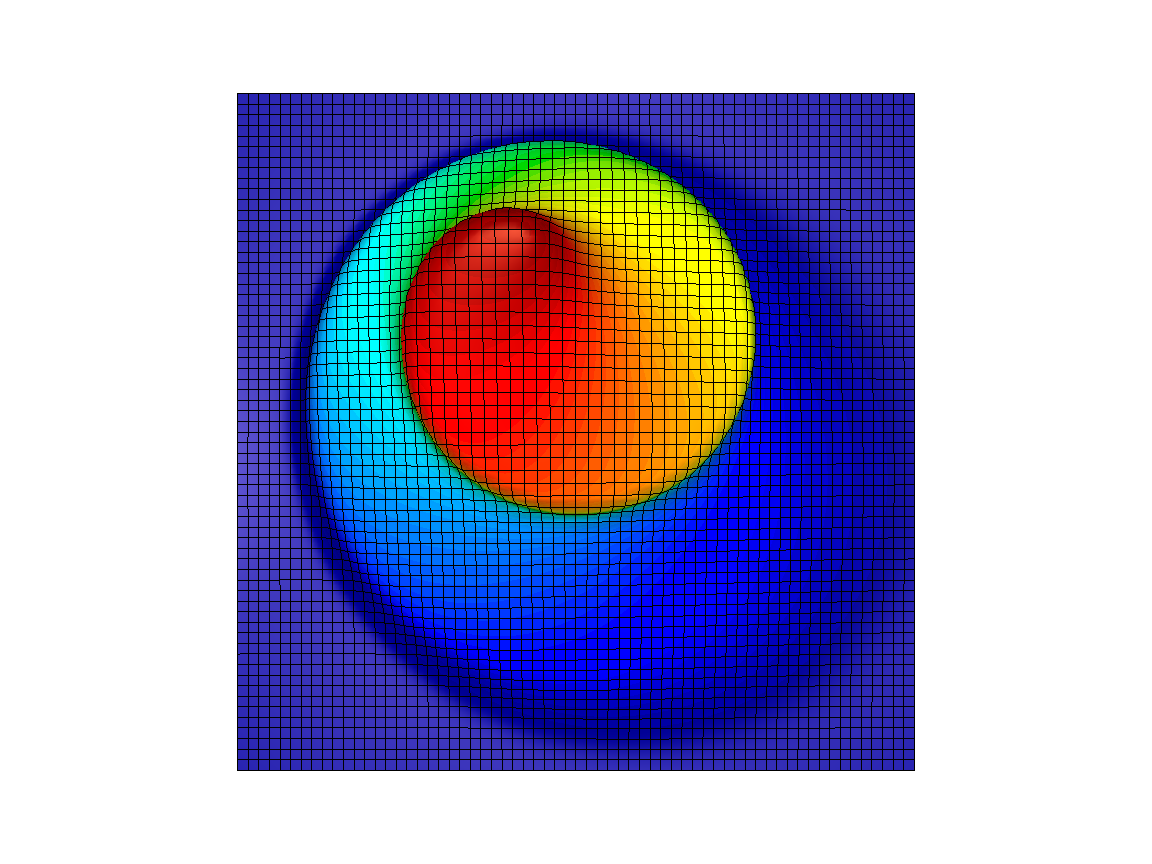}
		\caption{LO, $u_h \in [0.785,10.885]$}
	\end{subfigure}
	\begin{subfigure}[b]{.45\linewidth}
		\includegraphics[width=\linewidth,trim=0 50 0 0,clip]{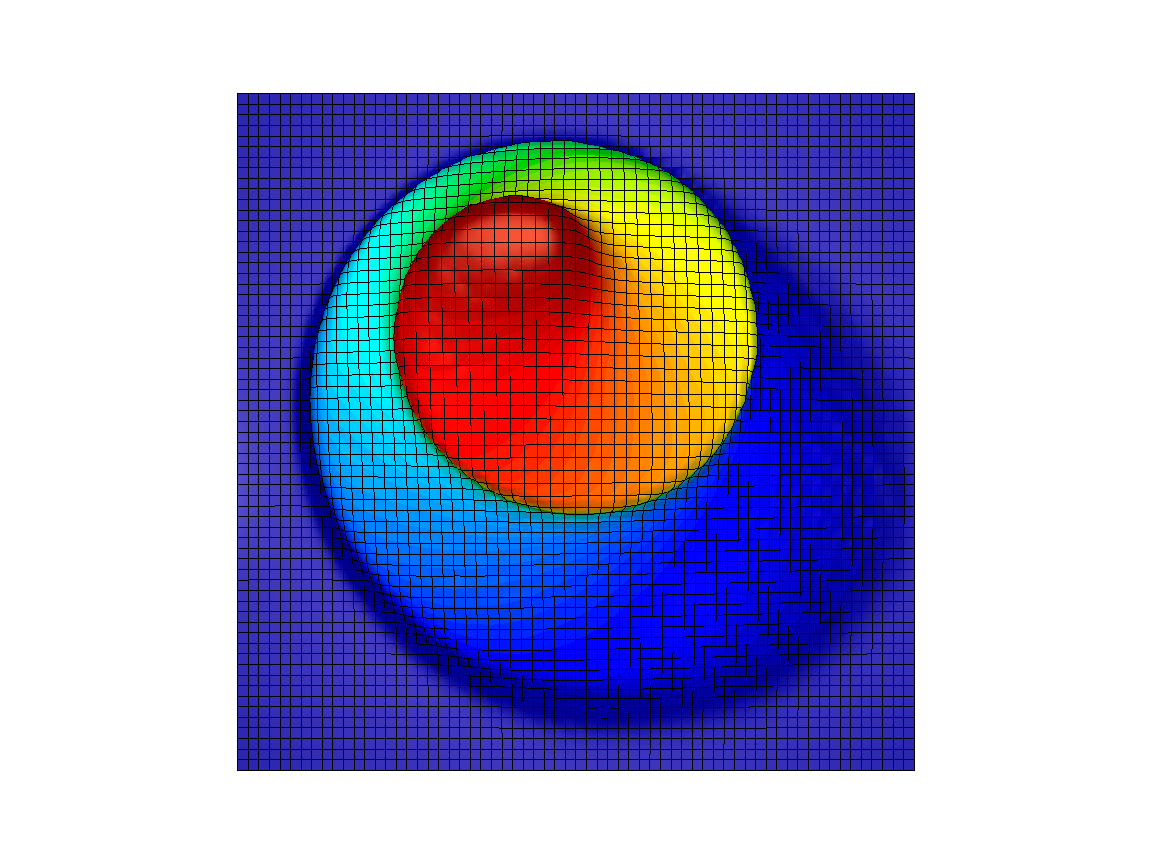}
		\caption{WENO, $u_h \in [0.785, 10.983]$}
	\end{subfigure}
	\caption{KPP problem, numerical solutions at $t=1$ obtained using $N_h=129^2$ and $p=2$.}
	\label{fig:kppp2}
\end{figure}

\begin{figure}[!htb]
	\centering
	\begin{subfigure}[b]{.45\linewidth}
		\includegraphics[width=\linewidth,trim=0 50 0 0,clip]{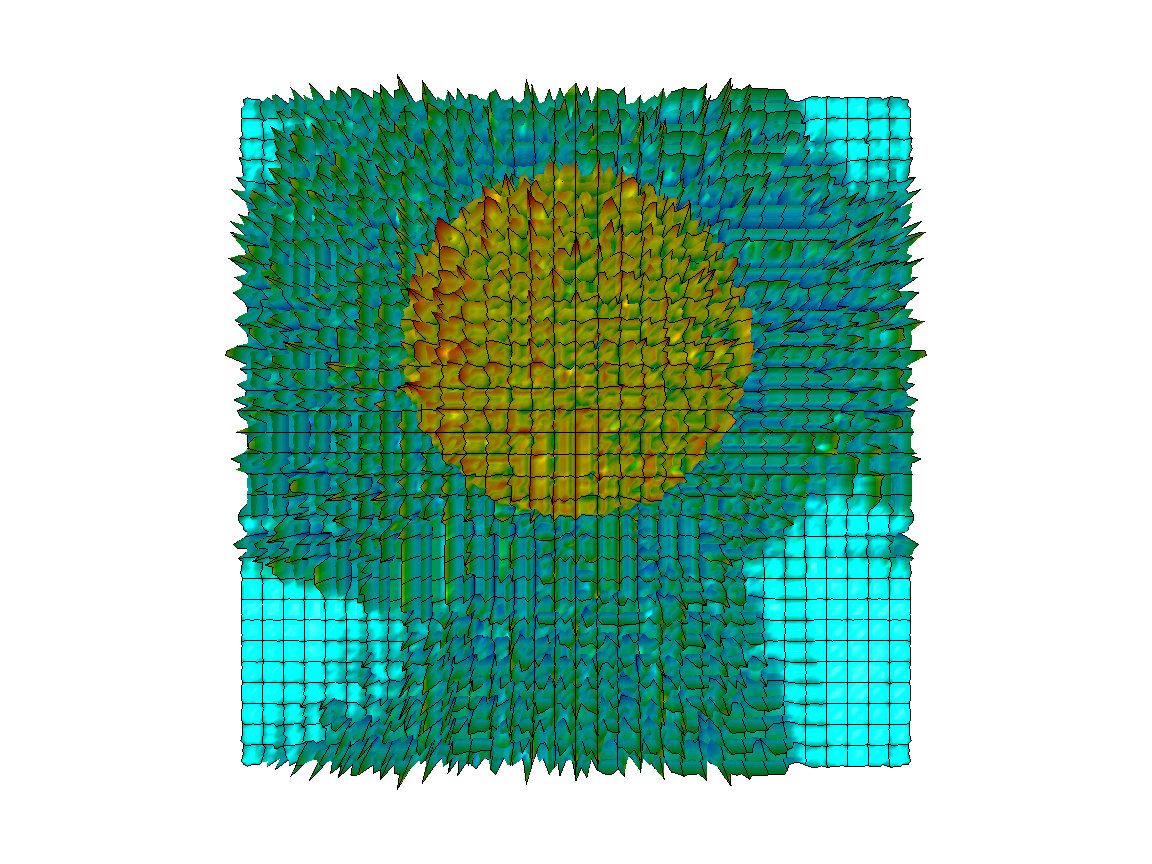}
		\caption{CG, $u_h \in [-16.260,28.369]$}
	\end{subfigure}
	\begin{subfigure}[b]{.45\linewidth}
		\includegraphics[width=\linewidth,trim=0 50 0 0,clip]{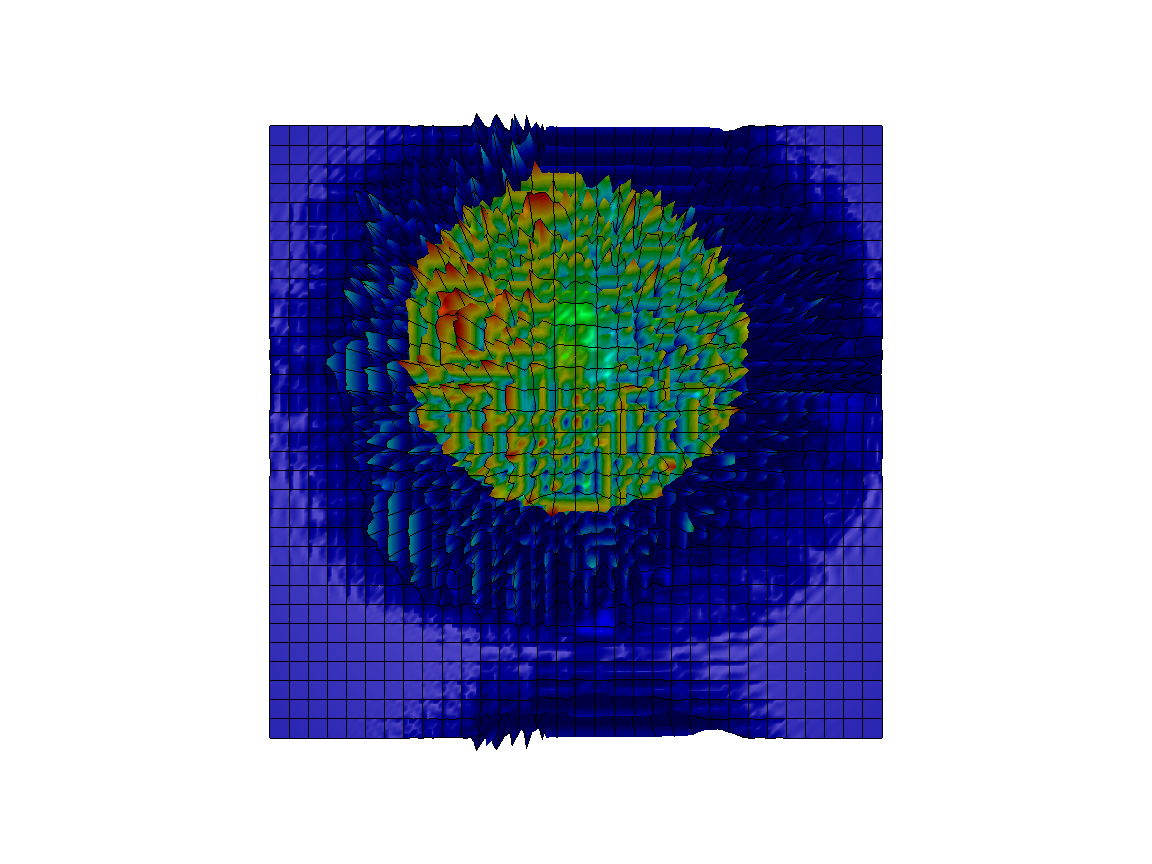}
		\caption{HO, $u_h \in [-10.828,22.094]$}
	\end{subfigure}
	\begin{subfigure}[b]{.45\linewidth}
		\includegraphics[width=\linewidth,trim=0 50 0 0,clip]{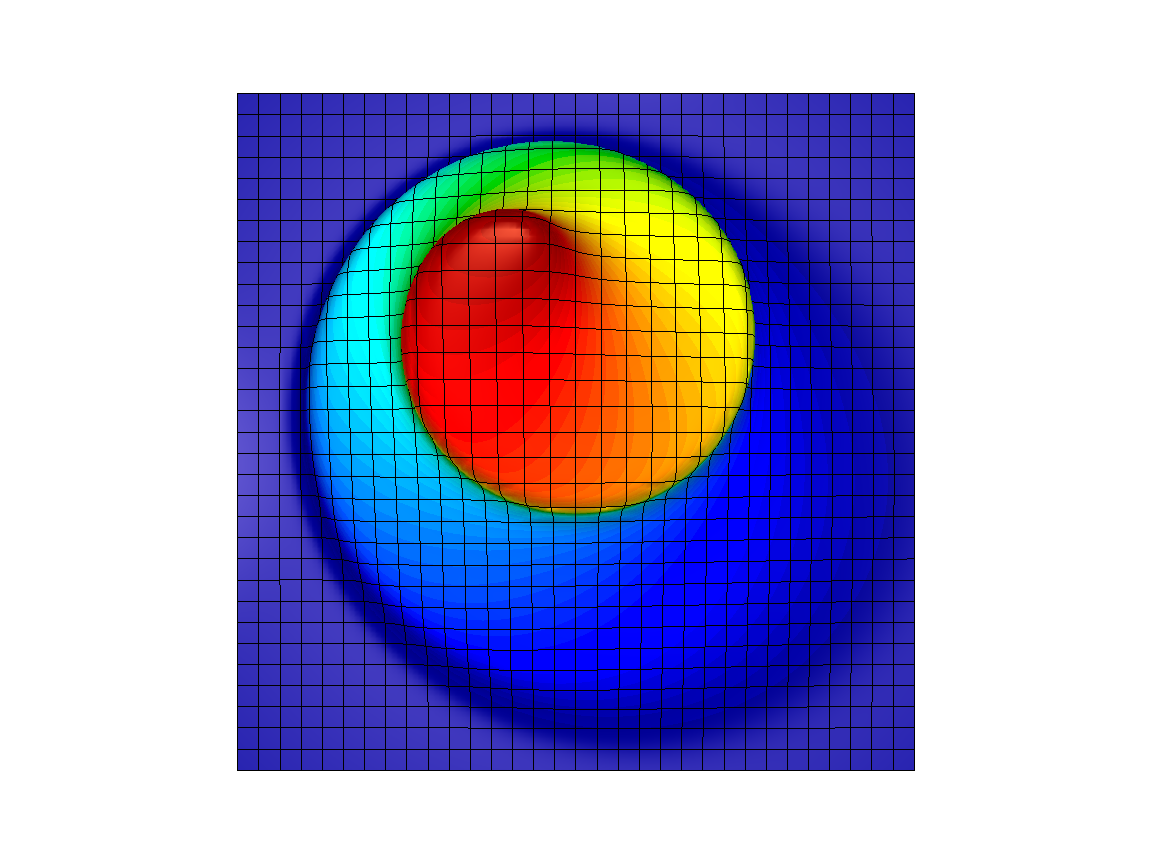}
		\caption{LO, $u_h \in [0.785, 10.882]$}
	\end{subfigure}
	\begin{subfigure}[b]{.45\linewidth}
		\includegraphics[width=\linewidth,trim=0 50 0 0,clip]{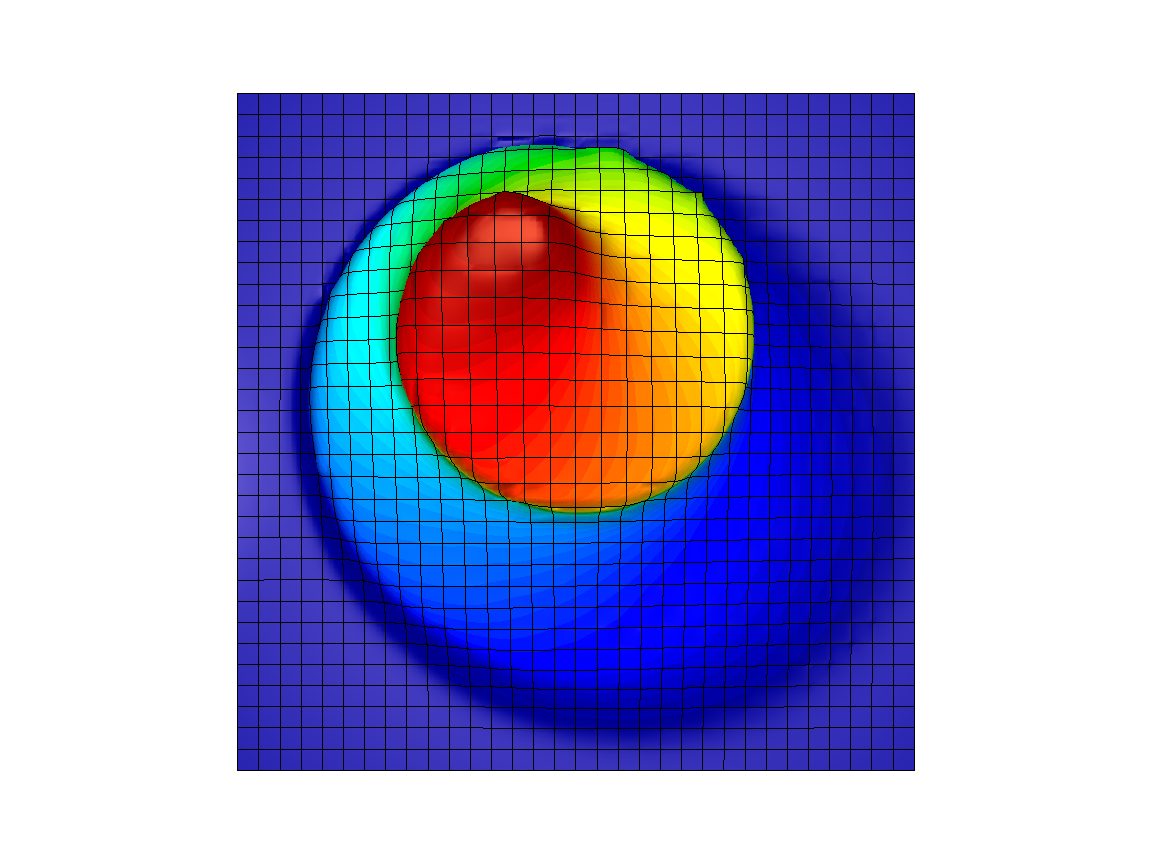}
		\caption{WENO, $u_h \in [0.778,10.980]$}
	\end{subfigure}
	\caption{KPP problem, numerical solutions at $t=1$ obtained using $N_h=129^2$ and $p=4$.}
	\label{fig:kppp4}
\end{figure}

As in the previous example,
we run simulations on uniform meshes of square cells corresponding to $N_h=129^2$ DoFs.
The numerical results for $p=\{1,2,4\}$ are displayed in Figs~\mbox{\ref{fig:kppp1}--\ref{fig:kppp4}}.
The lumped-mass CG method produces strong oscillations and an entropy-violating shock. The
HO stabilization alleviates the former problem (to some extent) but not the lack of entropy
stability. The LO results and the less diffusive WENO solutions are nonoscillatory. Moreover,
they reproduce the spiral wave structure of the entropy solution to the two-dimensional Riemann
problem. The small undershoot in the WENO result for $p=4$ is due to the use of a cell-based smoothness indicator on a coarse mesh. Such undershoots/overshoots disappear as the mesh is refined.
If necessary, preservation of global bounds and validity of entropy inequalities
can be strictly enforced using the algebraic flux correction tools presented in
\cite{KuHaRu2021,kuzmin2020g}. In fact, the development of our WENO scheme was largely motivated
by the possibility of using it as a baseline discretization for flux limiting approaches
that perform accuracy-preserving local fixes to enforce physical admissibility conditions
(cf.~\cite{kuzmin2023,zhang2010b,zhang2011}).



\section{Conclusions}
\label{sec:end}

We discussed a new way to combine high- and low-order components
of dissipation-based stabilization operators for CG discretizations of hyperbolic
conservation laws. Using piecewise-constant blending functions that measure deviations from
a WENO reconstruction, we managed to achieve optimal convergence rates
while avoiding spurious oscillations both globally and locally.
The resulting hybrid
scheme has the structure of a nonlinear artificial diffusion
method. Therefore, it is easier to analyze than DG-WENO schemes in which
partial derivatives of approximate solutions are manipulated
directly. Our theoretical investigations provide an optimal
estimate for the consistency error of the nonlinear stabilization
and worst-case \emph{a priori} error estimates for steady
advection-reaction equations.
The modular design makes it easy to upgrade individual
building blocks (stabilization operators, smoothness indicators,
reconstruction procedures) of the presented algorithm step-by-step.
In particular, it is worthwhile to investigate if the accuracy
of CG-WENO schemes can be improved by using Lagrange interpolation
polynomials (as in \cite{luo2007}), additional stencils for Hermite
interpolation, and/or alternative definitions of the linear weights.
Another promising avenue for further research is the use of
dissipation-based WENO stabilization in the DG context. Last
but not least, the cost of calculating smoothness indicators 
needs to be reduced, e.g., by using troubled cell detectors
as in \cite{zhong2013,zhu2016}.

\medskip
\paragraph{\bf Acknowledgments}
This article is dedicated to the memory of Prof.~Roland Glowinski, a great
mathematician who provided guidance and support to the first author many times over a span of two decades. Roland's wisdom and kindness will always be remembered.

The development of the proposed methodology was sponsored by the
German Research Association (Deutsche Forschungsgemeinschaft, DFG) under grant KU 1530/23-3.


\end{document}